\documentclass[11pt]{article}
\usepackage [latin1]{inputenc}
\usepackage{pdfpages}
\usepackage{amsfonts}
\usepackage{amssymb}
\usepackage{amsmath}
\usepackage{amsbsy}
\usepackage{mathtools}
\usepackage{mathrsfs}
\usepackage{makeidx}
\usepackage{MnSymbol}
\usepackage{pb-diagram}
\usepackage{tikz-cd}
\newtheorem{thm}{Theorem}[subsection]
\newtheorem{pro}[thm]{Proposition}
\newtheorem{fact}[thm]{Fact}

\newtheorem{lem}[thm]{Lemma}
\newtheorem{defi}[thm]{Definition}
\newtheorem{cor}[thm]{Corollary}
            
\newtheorem{exa}[thm]{Example}
\newtheorem{nota}[thm]{Notations}
\newtheorem{definota}[thm]{Definition and Notation}

\def\scrC{{\mathscr{C}}}

\def\scrK{{\mathscr {K}}}

\def\scrS{{\mathscr {S}}}

\def\scrCa{\mathscr{C}}
\def\calA{{\cal A}}

\def\calE{{\cal E}}
\def\calF{{\cal F}}

\def\calI{{\cal I}}

\def\calK{{\cal K}}

\def\calM{{\cal M}}
\def\calP{{\cal P}}
\def\calS{{\cal S}}

\def\calX{{\cal X}}
\def\calY{{\cal Y}}
\def\calZ{{\cal Z}}
\def\calSIM{{\calS\calI\!\calM }}

\def\bfa{\mathbf{a}}
\def\bfb{\mathbf{b}}
\def\bfc{\mathbf{c}}
\def\bfd{\mathbf{d}}

\def\bfx{\mathbf{x}}

\def\b1{{\mathbf{1}}}
\def\bca{{\mathbf{ca}}}

\def\g{\noindent}

\def\p{\vskip5truept}
\def\pp{\vskip 10truept}
\def\pg{\p\g}

\def\lppar{\mathop{(\hskip-4pt(}}
\def\rppar{\mathop{)\hskip-4pt)}}

\newcommand{\ppar}[1]{\lppar\hskip-2pt #1\hskip-2pt\rppar}

\newcounter{example}
\renewcommand\theexample{\Alph{example}}
\Alph{example}
\newenvironment{example}[1][]{\refstepcounter{example}\par\medskip
\noindent\textbf{Example~\theexample. #1} \rmfamily}{\medskip}

\usepackage[colorlinks=true,
            linkcolor=blue,
            urlcolor=blue,
            citecolor=blue]{hyperref}

\usepackage{hyperref}

\begin{document}

\title{A category of the political
\p\p\p
\Large{Part I - Hom\'onoia}}
\author{Joseph Abdou\thanks{Centre d'Economie de la Sorbonne, UMR 8174, Universit\' e Paris 1 Panth\' eon Sorbonne,  CNRS  }}
\date{\today \footnote{This article is a  revised version of the  September 9, 2024 paper \emph{Toward a category of the political}  https://hal.science/hal-04860761 } }
\maketitle

%
%

%
%

\pg {\bf Abstract.}
This research aims at providing a mathematical model of the  organization  of the polity and its  transformation. For that purpose
we construct two categories named  respectively Political Configuration  and Political Foundation. Our construction depends on a couple of variables called the foundational pair. One variable, called the Base, consists of a finite number of members (agents), while the other, called the Ground, consists of a set of states that reflect all relevant interests/values/aspirations of the base members. An object of the Configuration, called p-formation, extends  the  notion of simplicial complex, and  a morphism, which expresses the recomposition of the base, extends  the  notion of simplicial map. An object of the Foundation, called  p-site, describes  the profile of the polity, that is, how the states of the ground are intertwined between the agents.  A morphism between political sites consists of a  pair of maps, namely  a Base map and a Ground map, satisfying appropriate conditions. Two functors relate the Foundation and  the Configuration: the Knit which attributes to each p-site a p-formation  and the Nerve which attributes to each p-site a simplicial complex.  In the opposite direction a functor, called Canon,  which attributes to any p-formation its canonical p-site, turns out to be in an  appropriate sense the inverse of the Knit and the Nerve.

\tableofcontents

\section{Introduction}

  This study aims at developing a model of the polity. Any theory of the political must, within a coherent conceptual framework, account for both the organization and the transformation of the political entity.   \p
 The polity is to be understood, in the first instance as an entity constituted  by a collectivity  of individuals who mutually recognize their membership and affirm their commitment  to coexistence within its confines.   Insofar as these members are, by definition, bearers of distincts and autonomous  wills, the possibility of coexistence presupposes the establishment of a common ground capable of reconciling  divergent  attributes and interests. The realization of such coexistence  represents  the  fundamental telos of political activity, while the vocation of political theory   lies in elucidating  the conditions under which such  order comes into being.  
   
With this approach in mind, a  mathematical model of conflict resolution and compromise was elaborated in \cite{AK2019}. In that framework, the polity was described in terms of its organization, namely as  a collection of coalitions -- referred to as the \emph{political structure} -- which was assumed to form a simplicial complex. A coalition that constitutes  a simplex  is interpreted as viable, in the sense that some potential agreement binds its members.   

Despite its simplicity the  model  permits the study  of several interesting  political notions such as  the existence of friendly delegations and more generally the propensity  to reach a compromise, expressed through successive rounds of friendly delegations that eventually culminate in a viable coalition. Whether static or dynamic any notion of political activity conceivable within this model must be endogenous, meaning that it depends solely on the structure of the simplicial complex.  Indeed the model  does not provide a theory of transformation of political structures. 

 A further step  was undertaken in  \cite{Mock2021}    where the topology of simplicial complexes  was  employed to yield  instructive insights into politically meaningful notions such as stability, mediation, and merging. Mediation, in particular, necessitates  a passage from  a  given political structure to another.  
 
These considerations indicate firstly that a model of the polity  must integrate within a single conceptual framework, both the organization and the transformation of political structures. A categorical framework is particularly well suited to capturing these two dimensions in a unified manner. Secondly,  the polity is not  reducible to the organisation of the individuals composing  the  collectivity, but an entity founded  on  two variables, the Base, i.e., the set of individual composing the collectivity,  and the Ground, i.e., the set of wills, aspirations, and interests relevant for the collectivity.  Therefore our construction  must involve two conceptually distinct levels.    
 
 The \emph{Configuration} level represents the composition of the collectivity in terms of  coalitions, and specifies  the  transformations  that such a composition may undergo.  Formally the Configuration category depends on one variable, the Base.  An  object of this category  called  a political formation, is a structure  over the Base that extends the  notion of simplicial complex, and  a morphism between two political formations is a transformation that  extends  the well-known notion of simplicial map.

 The \emph{Foundation} level represents  the profile of the collectivity, i.e., the distribution of wills, aspirations, and interests among  its members,  and specifies the transformations  that such a profile  may sustain. Formally the Foundation category depends on a pair of variables, the Base and the Ground. An object of this category called a political site, is a profile of the Base over the Ground, and a morphism between two political sites is a regulated transformation of the profile. 

Within this construction, the Configuration  of the polity is treated as an observable, while the role of the Foundation is to  explain  what is observed in the Configuration. This is achieved by studying the structural 
relationships between the two levels, which is precisely the function of functors.  
The \emph{Knit} is the principal functor from the Foundation to the Configuration; it sends a political site to a political formation, and maps  a morphism of the Foundation to  a morphism of the Configuration. 
Another association, the \emph{Nerve}, sends a political site  to its simplicial complex (the  ``nerve'' of the political site in the algebraic topology literature). A key finding of this study is the central role of  the Knit in structuring the whole framework, whereas the Nerve  turns out to be functorial  only  on  some restricted domains.  
In the opposite direction, the functor \emph{Canon}  is introduced; it associates to any political formation  its canonical political site.  Since the Knit of the canonical site of a formation is precisely this very formation, the functor Canon has the Knit as a left inverse.  Any transformation that occurs in  a political  site  is reflected in a recomposition of the political formation and conversely any observable recomposition of a political formation can be attributed to some transformation of some political site.

In this part of the study, there is an important restriction concerning admissible transformations and recompositions:  in formal terms our morphisms are  built on functions\footnote{Generalization to correspondences will be the object of the second part of the study}; this amounts to  admit only recompositions that represent union of parties of the Base  and modifications  that represent  merger of states of the Ground.  

The paper is organized as follows.
The first section, devoted to the objects of our categories, introduces the notions of p-site and p-formation. Two operations are defined on p-sites, the Knit  and the Nerve.  The Knit maps a p-site to a p-formation with the same base, the Nerve maps a p-site to a p-structure with the same base. In the opposite direction, the  canonical map associates to any p-formation a canonical p-site.  It turns out that the Knit  map establishes a bijection  between   canonical p-sites  and p-formations  the inverse of which is precisely the canonical map.  On the other hand  the Nerve map  establishes a bijection  between the set of perfect canonical p-site and the set of p-structures (simplicial complexes). In section 2 we describe the  P-configuration  level  by introducing the morphisms that map a p-formation to another. The category  $\mathbf{P}^2_\bullet$ (complexes) whose objects  are p-formations, and its subcategory   $\mathbf{SIM}_\bullet$ (simplicial complexes) whose objects  are p-structures  are formally defined. 
In section 3  we start the description of the P-foundation level with either invariant Base or invariant Ground, by introducing the morphisms that map a p-site to another. In fact we introduce two classes of transformations, namely B-maps and S-maps, and we define the notions of B-isomorphism and S-isomorphism between p-sites correspondingly.  In section 4  we  study the general case with variable Base and Ground. We define the categories $\mathbf{bpol}_\bullet$ and $\mathbf{spol}_\bullet$,  the functors Knit, Canon, and Nerve, and establish the structural relationships between the two  levels precisely between 
 $\mathbf{bpol}_\bullet$ and p-configuration category $\mathbf{P}^2_\bullet$  on one hand, and $\mathbf{spol}_\bullet$ and p-configuration category $\mathbf{SIM}_\bullet$. 
In section 5 we focus on the scope of our model for applications and discuss some examples. 
In section 6 we conclude this part of the study by showing the need for an extension, thus motivating the production of the subsequent part.

\section{Political formation and  Political site}
Our concept of the polity depends on two variables , the Base denoted $I$, and the Ground denoted $A$. The couple $(I, A)$ will be referred to the as the \emph{foundational pair} of the polity. 
 The variable Base  $I$ is the set of individuals that compose the collectivity. Other denominations  for the Base would be the Board, the Constituency, the Collectivity, the Community. The institution of the polity determines who counts as member of the collectivity, and those who count form the Base.   An element of the Base will be called an \emph{agent}. Other denominations would be :  member, citizen, party.  The institution of the polity also determines  the physical, symbolic, and imaginary  values, that are acknowledged for members of the Base.   We call  \emph{Ground}  the set  of acknowledged values.  An element of the Ground will be called a \emph{state}.
 
 Objects of the Political Configuration  are political formations (p-formations, in the sequel), and depend  on the Base $I$, while  objects of  the Political Foundation are  political sites  (p-sites, in the sequel), and they depend on the foundational pair $(I, A)$.

%
%
%
%
%
%
%
%
%
%
%

\subsection{Preliminaries}
In what follows we assume that the Base $I$ in the foundational pair $(I, A)$ is a finite set.
\begin{nota}\rm\pg
 (i) For any finite set $X$, the set of all subsets of $X$ will be denoted $\calP(X)$.  
 \pg (ii) We put $\calP_0(X):=\calP(X)\setminus \{\emptyset\}$.  If  $\calF\subset \calP(X)$,  we put $\calF_0:=\calF\setminus\{\emptyset\}$,  $\calF^\downarrow := \{ s\in \calP_0(X)\vert \exists m\in \calF, s\subset m\}$ and  we denote  by $\calF^{max}$  the set of maximal elements (for inclusion) of $\calF$.  
\pg  (iii) If $f: X\to Y$ is a map in sets, ie. a function of domain $X$ and codomain $Y$,  we denote by $\overline{f}$ the corestriction  of $f$, that is,   $\overline{f}: X \to f(X)$ such that $\overline f(x) = f(x)$ $(x\in X)$, 

\pg (iv) If $\varphi : I \to J$  is a function, we denote by $\hat \varphi $ the  natural extension of $\varphi$ to subsets of $I$, precisely $\hat \varphi : \calP(I) \to \calP(J)$ defined by : $\hat \varphi (s) =\{\varphi (i) : i\in s\}$
\end{nota}

\subsubsection{Complexes and simplicial complexes}
  By  \emph{$I$-complex} (or a \emph{complex} on $I$) we mean  any subset  $\calK$ of $\calP(I)$. 
An element $s$ of $\calK$  is called a \emph{coalition}.
The set of all $I$-complexes,   $\calP(\calP(I))$,   will be denoted $\calP^2(I)$.
\pg  An \emph{$I$-simplicial complex} (or a \emph{simplicial complex} on $I$)  is  an $I$-complex  $\calK$  such that  $\calK^\downarrow=\calK$.  This is equivalent to say that $\emptyset \notin \calK$ and $ (s\in \calK, s'\subset s \Rightarrow s'\in \calK)$. Note that for any $I$-complex $\calK$, $\calK^\downarrow$ is the smallest simplicial complex on $I$ that contains $\calK_0$, called the simplicial complex generated by $\calK$. The set of all $I$-simplicial complexes is denoted  $\calSIM(I)$.
\pg The \emph{carrier }of a $I$-complex $\calK$ is the union of all coalitions in $\calK$

\p
\g Note that our definition of simplicial complex does not require that all singletons of $I$ be coalitions, in other terms we do not require that the carrier of $\calK$ be equal  to $I$, as is the case in the classical definition.
Note also that a coalition of a simplicial complex is classically called a \emph{simplex}.

\subsubsection{Profiles and Parting functions}

 \pg A \emph{profile }of $(I, A)$  is a function $\alpha: I\to \calP(A)$.  
The set of all profiles on $(I,A)$, $\calP(A)  ^I$, is denoted $\mathrm{PROF}(I,A)$. 
\pg A  \emph{parting map}  of $(I,A)$ is  a function $\rho: A\to \calP(I)$. The set of all parting maps, $\calP(I)^A$ is denoted $\mathrm{PART}(I,A)$. 

\pg  The parting map associated to the profile $\alpha$, is the function  $\Pi^{I,A}_{\alpha}: A \to \calP(I)$ defined by: 

\begin{equation}\label{parting1} \Pi^{I,A}_{\alpha}(a) =\{i\in I\vert a \in \alpha(i)\} \, (a\in A)
\end{equation} 
\pg The  profile associated to the parting map $\rho$,  is the function $\Gamma^{I,A}_{\rho}: I\to \calP(A)$  defined by : 
\begin{equation}\label{profiling}
\Gamma^{I,A}_{\rho} (i)=\{a\in A \vert\, i\in \rho(a)\}\, (i\in I)
\end{equation}

\pg Those associations determine two maps:

\begin{equation}
\begin{tikzcd}
{\mathrm{PROF}(I,A)}\arrow[r, shift right=0.5ex,  "\Pi^{I,A}" ']&{\mathrm{PART}(I,A)}\arrow[l,shift right=0.5ex, "\Gamma^{I,A}" ']
\end{tikzcd}
\end{equation}

\begin{fact} For any $(I,A)$:
\pg (i) $\Gamma^{I,A}= \Pi^{A,I}$
\pg (i) $\Pi^{I,A}$ and $\Gamma^{I,A}$ are inverse of each other. 
\end{fact}
\pg \emph{Proof.} $\rho=\Pi_\alpha \Leftrightarrow  \forall a \in A, \rho(a)= \{i \in I \vert\, a\in \alpha(i)\}$
\pg  $\Leftrightarrow \forall a\in A, \forall i\in I, ( i\in \rho(a)) \Leftrightarrow (a\in \alpha(i))$
\pg  $\Leftrightarrow \forall i\in I, \alpha(i)=\{a\in A \vert\, i\in \rho(a)\}$ $\Leftrightarrow \alpha=\Gamma_{\!\rho}$ \hfill $\Box$

 \subsection{P-formation}
Given a polity over some base,  a p-formation consists of a collection  
 of coalitions of members that share exclusively some well-identified interest. If for instance, the base is composed of citizens, the coalitions  are interest groups, such as guilds, clubs, corporations, trade unions, parties. If the base is composed of oligarchs then the coalitions are alliances. If the base is composed of parties as in representative democracies, then such coalitions are the potential coalitions that may form the government. 
  Note that the notion of p-formation may be considered as an extension of the notion of p-structure introduced in \cite{AK2019}.

\begin{defi}\label{defpformation}\rm
A \emph{p-formation} is a couple $(I,\calA)$ where $\calA\in \calP^2(I)$.  $I$ is the \emph{Base} and $\calA$ is the \emph{Complex}  of the p-formation. 

 \pg A \emph{p-structure}  $(I,\calE)$ is a p-formation where $\calE$ is an $I$-simplicial complex. \hfill$\triangledown$
\end{defi}

\g The notion of p-structure  is inherited  from  \cite{AK2019}.

\subsection{P-site}
\begin{defi}\label{defpsite}\rm
A  \emph{p-site} is a triple $\ppar{I, A, \alpha}$ where $I$ is a finite set, $(I,A)$ is a foundational pair and  $\alpha$ is a profile on $(I,A)$. \hfill $\triangledown$
\end{defi}

 \pg In the sequel we shall denote a p-site based on the pair $(I,A)$ by the boldface $\bfa$.
The profile $\alpha$ will be frequently denoted as a family $(A_i)$ where $A_i=\alpha(i)$ $(i\in I)$. Such a family is sometimes called  a cover of $A$. The notion of cover is classical in the algebraic topology literature  (e.g. \cite{Spanier} chapter 3, section I, example 6 and \cite{Barmak1} chapter 5, section 3).

\pg A p-site $\bfa\equiv \ppar{I, A, \alpha}\equiv\ppar{I, A, (A_i)}$ can also be characterized  by its parting map $\Pi^{I,A}_\alpha$  denoted $\pi_\bfa$. One has  $A_i =\{x\in A\vert \pi_\bfa(x)\ni i\}$, $\pi_\bfa(a)=\{j\in I\vert  A_j \ni a\}$.

\pg We denote by $\mathrm{PSITE}(I,A)$ the set of all p-sites that depend on $(I,A)$. It is clear that $\mathrm{PSITE}(I,A)$ is the same set  mathematical object as  $\mathrm{PROF}(I,A)$.

\begin{nota} \label{notations1} \rm 
 (i) If  $s \in \calP(I)$, $D\in  \calP(A)$ we use in the sequel  the notations:  $s^c := I \setminus S$,  $D^c= :A\setminus D$. If $(A_i)$ is the profile of $\bfa$ we put:  
  $A_s: = \cap_{i\in s} A_i$, 
  $ A^c_s: = \cap_{i\in s} A^c_i$, in particular $A_\emptyset= A^c_\emptyset =A$. We also put   $\overset{\frown}A_s : = A_s \cap A^c_{s^c}$, in particular $\bar A_\emptyset =A\setminus \cup_{i\in I} A_i$.
  \pg (ii) We adopt similar notations for p-sites $\bfb =\ppar{J, B, (B_i)}$, $\bfc\equiv\ppar{K,C,(C_k}$, $\bfx=\ppar{U, X, (X_u)}$ etc, by simply replacing letter $A $ by $B$,  $C$, $X$, letter $I$ by $J$, $K$, $U$ etc.
\end{nota}

\begin{fact} Let  $\bfa\equiv \ppar{I,A, (A_i)}$  be a p-site.  For any $s\in \calP(I)$ :
 \pg $\overset{\frown}A_s = \{ a \in A\vert\, \pi_\bfa(a)=s\}$, 
 $A_s =  \{ a\in A\vert\,  \pi_\bfa(a)\supset s\} = \sqcup\{ \overset{\frown}A_t:  t\supset s\}$.
\end{fact}

%
%

\pg \subsection{Knit and Nerve  of a p-site}
In this subsection we establish the main relationships between the p-sites (the objects of the P-foundation) and p-formations (the objects of the P-configuration).  In particular we show how any p-site with base $I$ induces a couple of interesting complexes on $I$, say $(\calS, \calK)$ such that $\calK^-=\calS$  and how conversely, to any such couple of complexes, we can associate, in a canonical way, a  p-site that in turn induces precisely those very complexes.  

\begin{defi}\label{defknitnerve}\rm
 To any p-site ${\bfa} =\ppar{I, A, (A_i)}$ we associate the following $I$-complexes:
 \pg (i) the \emph {Knit} of $\bfa$ :  
\pg $\kappa_\bfa:= \{s \in \calP(I) \vert \overset{\frown}A_s \neq \emptyset\}= \{s\in \calP(I)\vert\, \exists a\in A, \pi_\bfa(a)= s\}$. 

 \pg (ii) the \emph{Nerve} of $\bfa$ :
\pg  $\sigma_\bfa:= \{s \in \calP_0(I) \vert A_s\neq \emptyset\}=\{s\in \calP_0(I)\vert\, \exists a\in A, s\subset\pi_\bfa(a) \}$. \hfill$\triangledown$

\end{defi}

Note that the Nerve of $\bfa$ is a simplicial complex, therefore $(I, \sigma_\bfa)$ is a p-structure.  The notion of Nerve is classical in the algebraic topology literature (see e.g. \cite{Spanier} chapter 3, section I, example 6, and \cite{Barmak1} chapter 5, section 3).

\begin{fact}   For every p-site $\bfa$  we have   $ \sigma_\bfa   = (\kappa_ \bfa)^\downarrow$;  in other words $ \sigma_\bfa $ is the simplicial complex generated by $\kappa_ \bfa$.
  \end{fact}
It follows that  $\sigma_ \bfa$ is fully determined by $\kappa_\bfa$. On the other hand $\kappa_\bfa$ is not  determined by $\sigma_\bfa$.  
\begin{example} As an illustration of the Nerve and the Knit we present two p-sites on the same ground of 6 elements.  Both have the same Nerve, namely the triangle, a simplicial complex of dimension $1$,  composed of  all subsets of  cardinal 1 and 2. The Knit of the p-site on the left is equal to its Nerve, while the Knit of the p-site on the right is composed of subsets of cardinal 2.  
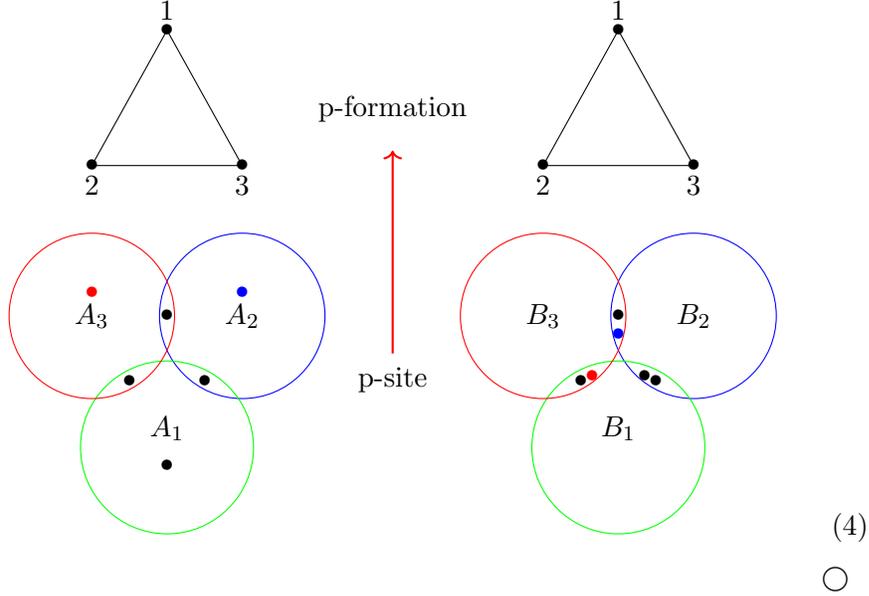
\begin{equation}
\begin{tikzpicture}[scale=0.5]

\draw (8,4) -- (10,7.6) -- (12,4) -- cycle ;
\draw (8,4) node{$\bullet$} node [below] {$2$};
\draw (10,7.6) node{$\bullet$}node [above] {$1$};
\draw (12, 4) node{$\bullet$} node [below] {$3$};
\draw  (20,4) -- (22,7.6) -- (24, 4) -- cycle ;
\draw (20, 4) node{$\bullet$} node [below] {$2$};
\draw (22,7.6) node{$\bullet$} node [above] {$1$};
\draw (24, 4) node{$\bullet$} node [below] {$3$};

\draw[->, red, thick] (16,-1) -- (16, 4.4);
\draw(16, -1.7) node{p-site};
\draw(16,5.5) node{p-formation};

\draw[blue] (12,0) circle (2.2);
\draw[green] (10,-3.5) circle (2.3);
\draw[red] (8,0) circle (2.2);
\draw(8,0) node{$ A_3$};
\draw(12,0) node{$ A_2$};
\draw(10,-3) node{$ A_1$};
\draw(10,0) node {$\bullet$};
\draw[red] (8, 0.60) node {$\bullet$};
\draw[blue] (12, 0.60) node {$\bullet$};
\draw(9,-1.75) node {$\bullet$};
\draw(11,-1.75) node {$\bullet$};
\draw(10,-4) node {$\bullet$};

\draw[blue] (24,0) circle (2.2);
\draw[green] (22,-3.5) circle (2.3);
\draw[red] (20,0) circle (2.2);
\draw(20,0) node{$ B_3$};
\draw(24,0) node{$ B_2$};
\draw(22,-3) node{$ B_1$};
\draw(22,0) node {$\bullet$};
\draw[red] (21.3,-1.60) node {$\bullet$};
\draw (22.7,-1.60) node {$\bullet$};
\draw[blue] (22, -0.50) node {$\bullet$};
\draw(21,-1.75) node {$\bullet$};
\draw(23,-1.75) node {$\bullet$};
\end{tikzpicture}
\end{equation}
\hfill $\bigcirc$
\end{example}
\begin{example}[]  \label{gallopolis}{(Gallopolis)}
Political life in \textit{Gallopolis} is organized around three fundamental dimensions: the economic issue, the societal issue, and the openness issue.  

The \textbf{economic issue} concerns the level of taxation and the extent of public services. It is represented by the set  
\[
E = \{1,2,3\},
\]  
where $1$ denotes low taxation and limited public services, and higher values correspond to higher levels of taxation and public provision.  

The \textbf{societal issue} relates to matters such as family, abortion, and gender. It is represented by the set  
\[
S = \{l, n, c\},
\]  
where $l$ stands for \textit{liberal}, $n$ for \textit{neutral}, and $c$ for \textit{conservative}.  

The \textbf{openness issue} concerns immigration, attitudes toward foreigners, and xenophobia. It is represented by the set  
\[
O = \{\alpha, \beta, \gamma\},
\]  
where $\alpha$ denotes a policy of zero immigration and strictness toward immigrants, while $\beta$ and $\gamma$ correspond to progressively more open policies.  

The political space of Gallopolis is therefore defined as the product set  
\[
A = E \times S \times O.
\]  

Within this space, five political parties constitute the principal political forces, each characterized by a set of desirable states:
\begin{align*}
\text{LEFT}   &= \{2,3\} \times \{l\} \times \{\gamma\}, \\
\text{SOCD}   &= \{2\} \times \{l,n\} \times \{\beta,\gamma\}, \\
\text{CONS}   &= \{1,2\} \times \{n,c\} \times \{\alpha,\beta\}, \\
\text{LIBER}  &= \{1\} \times \{l\} \times \{\alpha,\beta,\gamma\}\cup \{1\} \times \{n\} \times \{\beta,\gamma\}\\
\text{RIGHT}  &= \{1,2\} \times \{c\} \times \{\alpha\}.
\end{align*}

\pg$\bullet$
The nerve of Gallopolis is the following simplicial complex represented by its geometric realization. 

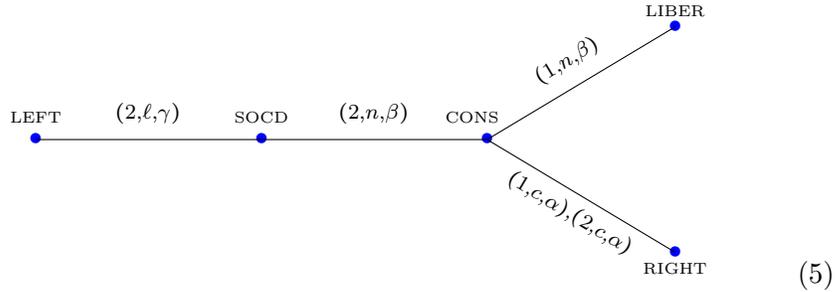
\begin{center}
\begin{equation}\label{figgallopolis1}
\begin{tikzpicture}[scale=1]
\draw[blue] (0,0) node{$\bullet$}; \draw (0,0.3)node  {\tiny LEFT};
\draw[blue] (3,0) node{$\bullet$}; \draw  (3, 0.3) node {\tiny SOCD};
\draw[blue] (6,0) node{$\bullet$}; \draw (5.8,0.3) node {\tiny CONS};
\draw[blue] (8.5,1.5) node{$\bullet$} node[above, black] {\tiny LIBER};
\draw[blue] (8.5,-1.5) node{$\bullet$} node[below, black] {\tiny RIGHT};
\draw(0,0)--(3,0)node[midway,above]{${}^{(2,\ell, \gamma)}$}--(6,0)node[midway,above]{${}^{(2, n, \beta)}$};
\draw(6,0)--(8.5,1.5)node[midway,above,sloped] {${}^{(1,n,\beta)}$};
\draw(6,0)--(8.5,-1.5)node[midway,below,sloped] {${}^{(1,c,\alpha), (2,c,\alpha)}$};
\end{tikzpicture}
\end{equation}
\end{center}
\pg
The nerve is a simplicial complex of dimension $1$.  In the figure  every simplex of dimension $1$ is affected by the states that are in the intersection of its vertices. 
\pp


\pg$\bullet$
The parting function  of Gallopolis $\pi\equiv \pi_\bfa$  can be partially read on graphic (\ref{figgallopolis1}),   precisely :

\pg - on the  states that connect two vertices we have:    $\pi (2,l,\gamma)= \{\mathrm{LEFT, SOCD}\}$, $\pi (2,n,\beta)= \{\mathrm{SOCD, CONS}\}$, etc.   

\pg - On the rest of the states the parting function is as follows:
\pg $\pi(3, l, \gamma) =\{\mathrm{LEFT}\}$,
\pg $\pi(1,l, \alpha) = \pi(1,l, \beta) =\pi(1,l, \gamma) =\pi (1,n, \gamma)= \{\mathrm{LIBER}\}$,
\pg $\pi(2, n, \gamma)=\pi(2, l, \beta)=\{\mathrm{SOCD}\}$,
\pg $\pi(1,c, \beta)=\pi(1,n, \alpha)=\pi(2,n, \alpha)=\pi(2,c,\beta)=\{\mathrm{CONS}\}$. 
\pg 
 -  the parting function takes the value $\emptyset$ on the states, that do not belong to any party, for instance $(1,c,\gamma), (3, l, \alpha), (3, l, \beta), (3, n, \alpha), \ldots $

\pg$\bullet$ The Knit of Gallopolis, $\kappa_\bfa$, is the image of the parting function  $\pi_\bfa$. It  has as maximal elements  precisely the simplexes of dimension $1$ of the nerve,  the other elements  being $\emptyset$, $\{\mathrm{LEFT}\}, \{\mathrm{LIBER}\}, \{\mathrm{SOCD}\},$ $ \{\mathrm{CONS}\}$. Therefore the Knit of Gallopolis  (including the empty set) is a complex  of cardinality $9$.  The only ``coalition" that is  present in the nerve but absent from the Knit is $\{\mathrm{RIGHT}\}$. 
 
 \pg Note that, in this model of the polity,  every party except RIGHT can be identified by characteristic  principles. For instance Left is the member of society that can be identified by its principle $(3,l,\gamma)$, LIBER  can be identified by a set of principles, namely $ \{(1,l, \alpha), \pi(1,l, \beta), \pi(1,l, \gamma), (1,n, \gamma)\}$; by contrast no principle can identify RIGHT. This property is a consequence of the fact that  the aspirations of RIGHT are a proper subset of the aspirations of CONS. \hfill $\bigcirc$ 

\end{example} 
\subsection{Some classes of p-sites}\label{familySC} 
\subsubsection{Canonical and subcanonical p-sites} 
Let $I$ be a base and  let $\calA \subset \calP(I)$ $\calA\neq \emptyset$.  $\calA$ will be considered as the ground of a family of p-sites called the subcanonical family of p-sites with base $I$.   For any  $\calF\subset \calP_0(I), \calF\subset \calA$, we define the p-site $\bfc^\calF(I,\calA)=\ppar{I, \calA, (\calA_i^\calF)}$ with base $I$, ground $\calA$ and with profile:
 \begin{equation} \calA^\calF_i: =\{ s\in \calF: s\ni i\} (i\in I) 
 \end{equation}
 
\pg Dually  $\bfc^\calF (I,\calA)$ can be defined as the  p-site $\bfa$ with ground $\calA$ and parting function:
 
 \[ \pi_\bfa(s) = \left\{
  \begin{array}{ll}
    s& : s\in \calF\\
    \emptyset & : s \in  \calA\setminus \calF
  \end{array}
\right.
\]
  $\bfc^F (I,\calA)$ can be characterized as the p-site  $\bfa$ with ground   $\calA$ such that:
 
  \[\bar X_s = \left\{
  \begin{array}{ll}
    \!\!\{s\}& : s\in \calF\\
    \emptyset & : s\in  \calP(I) \setminus \calF
  \end{array}
\right.
\]

\g Such a p-site is said to be \emph {subcanonical} with ground $\calA$ and \emph{effective ground} $\calF$.

\begin{pro} Let $\bfa := \bfc^\calF(I,\calA)$. Then:
\pg (i) $\kappa_{0\bfa}=\calF$.
\pg(ii) $\kappa_\bfa=\left\{
  \begin{array}{ll}
    \calF &  { \rm if\, } \calA = \calF \\
    \calF\cup\{\emptyset\} &{\rm  if \,}   \calA \neq  \calF
  \end{array}
\right.
$

\end{pro}

\begin{cor} Let $ \calF \subset \calP_0(I)$.  A subcanonical p-site $\bfa$ is such that  $\kappa_{0\bfa}=\calF$ if and only if  $\bfa= \bfc^{\calF} (I, \calA)$ for some $\calA$,   $\calA\supset \calF$.
\end{cor}

\begin{cor} Let $ \calK  \in  \calP^2(I)$.  What are the subcanonical p-sites $\bfa\equiv\bfc^\calF(I, \calA)$ such that $\kappa_\bfa=\calK$  ?
\pg  - If $\emptyset \in \calK$,  the solution is any one of the p-sites  where $\calF=\calK_0$ and $\calA\in \calP^2(I)$, $\calA\supset \calK_0$, $\calA\neq \calK_0$, in particular  the solution $\calF=\calK_0$, $\calA=\calK$. 
\pg - If  $\emptyset \notin  \calK$,   the  unique solution is the subcanonical p-site  where $\calF=\calA=\calK$.
\end{cor}

\begin{defi}\label{canonicaldefi}\rm
The \emph{canonical}  p-site   associated to the p-formation $(I, \calA)$ is  the p-site with foundational pair $(I, \calA)$ and profile $\calA_i=\{ s \in \calA\vert\, i\in s\}$.  It is denoted  $\bca(I,\calA)$.   A p-site $\bfa$ with base $I$ is said to be canonical on $I$, if $\bfa =\bca(I, \calA)$ for some $I$- complex $\calA$.  Clearly $\bca(I,\calA)\equiv \bfc^{\calA_0}(I,\calA)$.
 \end{defi}

 \begin{pro}\label{can} Let $\bfa =\ppar{I, \calA, (\calA_i)}$, where $\calA\subset \calP(I)$. The following are equivalent:
 \pg (i) $\bfa$ is canonical,
 \pg (ii)  $ \forall  s\in \calA, \pi_\bfa (s)=s$, equivalently $\pi_\bfa$ is the inclusion $i^\calA: \calA\hookrightarrow \calP(I)$,
 \pg (iii) $\overset{\frown} \calA_s =\{s\}$ for all $s\in \calA$, and $\bar \calA_s= \emptyset$ for all $s\in \calP(I)\setminus \calA$.
 \end{pro}

\begin{cor}  \label{canonical2}Let $\calK\in \calP^2(I)$. There is a unique canonical p-site $\bfa$ such that $\kappa_\bfa=\calK$, precisely $\bca(I,\calK)$.
\end{cor}
%

\begin{cor} \label{canonical3} Let $\calE\in \calS \calI \!\calM (I)$. The canonical p-sites $\bfa \equiv \bca(I, \calA)$ such that $\sigma_\bfa =\calE$ are those where $\calE^{max}\subset \calA\subset \calE$
\end{cor}

\begin{defi}  \rm A p-site $\bfa$ is \emph{perfect} if $\kappa_{\bfa} =\sigma_\bfa$ \hfill$\triangledown$
\end{defi}

\begin{cor} \label{canonical3} Let $\calE\in \calS \calI \!\calM (I)$.  There exists one and only one perfect canonical p-site  $\bfa$ such that $\sigma_\bfa= \calE$,  namely  $\bca(I, \calE)$%
\end{cor}

\subsubsection{Simple p-sites} 
\begin{defi}\rm  A p-site $\bfa=\ppar{I, A, (A_i)}$ is  said to be \emph{simple} if its parting function $\pi_\bfa : A \to \calP(I)$ is injective.
\end{defi}

\begin{defi}\rm  Let $f: A \to A'$ be a function, and let $\bfa=\ppar{I, A, (A_i)}$ and $\bfa'=\ppar{I, A', (A'_i)}$ be two p-sites with the same base $I$. $f$ is said to be an \emph{isotopy} of $\bfa$ to $\bfa'$ if $f$ is a bijection  and $f(A_i)= A'_i$ $(i\in I)$. 
If $f$ is an isotopy of $\bfa$ to $\bfa'$ then its inverse  $f^{-1}$ is an isotopy of $\bfa'$ to $\bfa$. Two p-sites $\bfa$ and $\bfa'$ are said to be \emph{isotopic} if  there is  an isotopy  between the two p-sites.
\end{defi}

\begin{fact} Any canonical p-site is simple. Conversely any simple p-site  $\bfa$ with base $I$ is isotopic to $\bca(I, \kappa_\bfa)$. In fact $\overline\pi_\bfa$ is an isotopy of $\bfa$ to $\bca(I, \kappa)$.
\end{fact}
\emph{Proof}. In view of  \ref{can}, the parting function of $\bca(I, \kappa_\bfa)$ is injective, therefore the p-site $\bca(I, \kappa_\bfa)$ is simple.  Conversely for any  simple p-site $\bfa = \ppar{I, A, (A_i)}$, the parting function $\pi_\bfa$ is a bijection $A \to \kappa_\bfa$. Put $\bfb:= \bca(I, \kappa_\bfa)$ and $\kappa_\bfa= \calA$. $\overline{\pi}_\bfa (A_i)= \pi_\bfa (A_i)=\{\pi_\bfa (x) : x \in A_i\} =  \{\pi_\bfa (x) : i \in \pi_\bfa (x)\} =\{s \in \kappa_\bfa \vert\, i\in s\} = \calA_i$. \hfill $\Box$

\begin{pro} For any finite set $I$, the Knit map $\kappa$ is a bijection between the set of all canonical p-sites over base $I$  and the set of all p-formations  over base $I$, $\calP^2(I)$, the inverse of which is the canonical map $\bca (I, \cdot)$. In symbols:

$\begin{tikzcd}
{\textrm{Canonical p-sites over base $I$} }\arrow[r, "{\kappa} ", yshift =0.7ex]& {\calP^2(I)} \arrow[l, " \bca"  , yshift=-0.7ex]
\end{tikzcd}$

\pg The Nerve  map is a bijection between the set of perfect canonical p-sites and the set of simplicial complexes on $I$, $\calSIM(I)$, the inverse of which is the canonical map
 $\bca (I, \cdot)$. In symbols:

$\begin{tikzcd}
{\textrm{Perfect canonical p-sites over base $I$} } \arrow[r, "{\sigma=\kappa} ", yshift =0.7ex]&
 {\calSIM(I)} \arrow[l, " \bca"  , yshift=-0.7ex]
\end{tikzcd}$

\pg \emph{Proof.} Straightforward\hfill$\Box$

\end{pro}

\section{P-configuration}
The p-configuration is the framework of the study of p-formations and their transformations. The transformations of p-formations are called \emph{recompositions}. In this part of our study a recomposition is built   on (we will also say ``induced by") a function  between the bases of the p-formations. By this construction the  p-configuration is established as a Category. As an application we show  in  Appendix \ref{delegation} how the notion of delegation, defined in \cite{AK2019} can be viewed as a  map in this category.

 \subsection{Recomposition of p-formations: BC-maps, SC-maps}
\begin{defi}\label{defmap2} \rm  Let $(I,\calX)$ and $(J,\calY)$ be two p-formations and let $\varphi: I\to J$ be a function. $\varphi$ is said to be a  \emph {BC-map} (resp. \emph{SC-map}) of domain $(I,\calX)$ and codomain $(J,\calY)$ if $\hat\varphi (\calX)= \calY$ (resp. $\hat\varphi (\calX)\subset \calY$)). \hfill$\triangledown$

\end{defi}

\begin{defi}\rm  Let $\varphi: I\to J$ be a function.  The  \emph{C- image}  of $(I,\calX)$ by $\varphi$ is the complex $(J, \hat \varphi (X))$. This is the 
 unique complex  $(J,\calY)$  that makes $\varphi : (I,\calX) \to (J,\calY)$ a BC-map.  This is also the smallest complex  $(J,\calY)$   for which  $\varphi$  is an SC-map from $(I,\calX)$ to $(J,\calY)$. \hfill$\triangledown$

\end{defi}

\begin{fact} \label{Ccomposition} Let $\varphi : (I,\calX) \to (J,\calY)$ and $\psi: (J,\calY) \to (K,\calZ)$ be BC-maps (resp. SC-maps)  then $\psi \circ \varphi$ is a BC-map  (resp. SC-map) with domain $(I,\calX)$ and codomain $(K,\calZ)$.
\end{fact}
 \emph{Proof.} Straightforward\hfill$\Box$


\pp\g
The language of category theory needed in this study can be found in any basic book on category theory (e.g. \cite{Leinster}, \cite{Awodey}  or \cite{Maclane}). 
\subsection{Categories of P-configuration:  $\mathbf{P}^2_\bullet$  and  $\mathbf{SIM}_\bullet$} \label{complexes}
 We define the category  $\mathbf{P}^2_\bullet$ (resp. $\mathbf{P}^2_=$)  whose objects  are p-formations.  Precisely an object  is an ordered pair  $(I,\calX)$, where $I$ is a finite set, and where  $\calX\in \calP^2(I)$. 
If $(I,\calX)$ and 
$(J,\calY)$  are two objects, a morphism   in  $\mathbf{P}^2_\bullet$ (resp. $\mathbf{P}^2_=$ ) of domain $(I,\calX)$ and codomain $(J,\calY)$ is an SC-map  $\varphi: I\to J$,  i.e., such that $\hat \varphi(\calX) \subset \calY$)  (resp. a  BC-map,  i.e., such that $\hat \varphi(\calX) = \calY)$.  The identity on $(I,\calX)$ is the one induced by $id_I$. In view of Fact \ref{Ccomposition}, the axioms of  category  (semi group law on morphisms, and  existence of identity)  are satisfied.

\pg Let $ \mathbf{SIM}_\bullet$ be the subcategory of $\mathbf{P}^2_\bullet$  with objects consisting only of  p-structures $(I,\calE)$ (simplicial complexes) and  whose  morphisms  are those inherited from $\mathbf{P}^2_\bullet$.  Note that   $ \mathbf{SIM}_\bullet$ is the usual category of simplicial complexes with simplicial maps as morphisms. (see \cite{Spanier} chap. 3)

\subsubsection{Subcategories  $\mathbf{P}^2_\subset[I]$, $\mathbf{SIM}_\subset[I]$.}
Let  $(\calP^2(I), \subset)$ be  the set of all subsets of $I$-complexes partially ordered by inclusion. We view this poset  as a category denoted  $\mathrm{P} ^2_\subset(I)$ whose objects are $I$-complexes, i.e. elements of $\calP^2(I)$, denoted   $\calX$, $\calY$, ... , and where there is a unique morphism, represented by the inclusion  $i^\calX_\calY: \calX \hookrightarrow \calY$, if  $\calX$ is included in $\calY$ and no morphism if $ \calX$  is not included in $\calY$.  
\pg
$\mathbf{P} ^2_\subset(I)$  is a subcategory of $\mathbf{P} ^2_0$ whose objects are $I$-complexes, and where there is a unique morphism  of domain $\calX$ and codomain $\calY$ if $\calX\subset \calY$, namely the unique SC-map induced by the identity map on $I$, and where  there is no morphism if $\calX\nsubset \calY$.

\pg
 We denote by $\mathbf{SIM}_\subset[I]$ the subcategory of  $\mathbf{P}^2_\subset [I]$  with objects the $I$-simplicial complexes.

\section{P-foundation}
The objects of the  P-foundation category are p-sites.  Maps between p-sites $\bfa\equiv\ppar{I,A, (A_i)}$ and $\bfb \equiv\ppar{J,B,(B_j)}$ are based on functions between foundational pairs $(I,A)$, and $(J,B)$, namely $\varphi: I\to J$ and $f:A\to B$. In order to understand how  the Nerve and the Knit  behave in such transformations, we are led to  define two types of transformations (and therefore two different categories) corresponding either to a conservative transformation (B-map\footnote{prefix B in reference to Boole}) or to  an expansive transformation (S-map\footnote{prefix S in reference to Simplex}). This construction leads to the definition of the two categories of P-foundation $\mathbf{bpol}_\bullet$ and  $\mathbf{spol}_\bullet$ and to the Functors Nerve and Knit that relate the P-foundation to the P-configuration.

For the clarity of the exposition, we study first the case where the Ground is invariant, i.e., $A=B$, $f= id_A$, then the case where the Base is invariant $I=J$ and $\varphi=id_I$, and finally the case where both the Base and the Ground are variable will be developed in the next section.

\subsection{P-foundation on invariant ground} \label{fixedground}

In this subsection we consider the class of   p-sites with ground  a set $A$, and we define transformations between two such p-sites that keep  the ground $A$ invariant.

\subsubsection{BP-maps, SP-maps}
Let $\varphi :I\to J$ be a function. In what follows $I$ and $J$ represent bases of p-sites $\bfa$ and $\bfb$ respectively.  Intuitively base $J$ succeeds to base $I$, $J$ being the result of the recomposition of $I$ through $\varphi$.  We interpret $f(i)=j$ as the fact that $i$ contributes to the emergence of $j$. If $f^{-1} (j) = \{i_1, \ldots,  i_p\}$ we interpret $j$ as a new political entity or party that is the result of the union or merger of  $i_1, \ldots,  i_p$, its parents. In order that this recomposition be simplicial the function $\varphi$ has to satisfy  structural conditions that make it a map of domain the p-site $\bfa$  and codomain the p-site $\bfb$. The following definition conveys through algebraic properties the admissible recompositions of the political landscape with an invariant ground. 

\begin{defi}\label{defmap}\rm  Let $\bfa\equiv \ppar{I, A, (A_i)}$ and $\bfb\equiv\ppar{J,A, (B_j)}$  be two p-sites \emph{with the same ground} $A$, and let  $\varphi: I \to J$ be a function. We define two classes of maps  with domain (or source)  $\bfa$ and codomain (or target) $\bfb$.
\pg $\varphi$ is a  \emph{BP-map}  if :

\begin{equation} \label{Bmap1}
\underset{i \in \varphi^{-1}(j)} \cup A_i  =  B_j\,\,\, ( j\in J)
\end{equation}

\pg $\varphi$ is a  \emph{SP-map} if :
\begin{equation}
\underset{i \in \varphi^{-1}(j)} \cup A_i  \subset  B_j\,\,\, ( j\in J)
\end{equation}
\end{defi}
\hfill$\triangledown$
\begin{nota}\rm If $\varphi: I\to J$ is any function, we denote by $\hat \varphi$ its natural extension to a map $\calP(I)\to \calP(J)$, precisely $\hat \varphi (S) =\{\varphi(i): i\in S\}$ $(S\in \calP(I))$.
\end{nota}
\begin{pro}\label{BGcharacter} $\varphi:I\to J$ is a BP-map (resp. an SP-map) if and only if $ \hat \varphi \circ \pi_{\bfa}= \pi_{\bfb} $ (resp. $ \hat \varphi \circ \pi_{\bfa}\subset \pi_{\bfb} $)
\end{pro}
\emph{Proof.}   We  give an extended proof for the B-case. Assume $ \hat \varphi \circ \pi_{\bfa}= \pi_{\bfb} $. Then $\forall j\in J, \forall a\in A$:
\pg
$a\in B_j \Leftrightarrow j\in  \pi_\bfb(a) \Leftrightarrow^*  j\in \hat\varphi(\pi_\bfa(a)\Leftrightarrow  \exists i\in \pi_\bfa(a), \varphi(i)=j  \Leftrightarrow \exists i\in I, a\in A_i, \varphi(i)=j \Leftrightarrow a\in \cup_{i\in \varphi^{-1}(j)} A_i$
\pg  This proves : $B_j =  \cup_{i\in \varphi^{-1}(j)} A_i$ for all $j$. It follows that  $\varphi$ is a BP-map.
\pg Conversely, assume that $\varphi$ is a BP-map: $\forall j\in J, \forall a\in A$:

\pg $ j\in \hat\varphi(\pi_\bfa(a) \Leftrightarrow \exists i\in \pi_\bfa(a), \varphi(i)=j \Leftrightarrow  j \in J, a\in \cup_{i\in \varphi^{-1}(j)} A_i
\Leftrightarrow^{**}  j \in J, a\in B_j \Leftrightarrow j\in \pi_\bfb(a)$
\pg This proves $\pi_\bfb(a) = \hat\varphi(\pi_\bfa(a))$ for all $a\in A$. It follows that  $\pi_\bfb = \hat\varphi\circ \pi_\bfa$
\pg A proof of the S-case is obtained by replacing  in the above lines $\Leftrightarrow^*$ by $\Leftarrow$, and the  $\Leftrightarrow^{**}$  by $\Rightarrow$. \hfill $\Box$

\begin{cor}\label{image2} Let $\varphi:I\to J$ a function and let $\bfa=\ppar{I, A, (A_i)}$, $\bfb=\ppar{J, A, (B_j)}$.
\pg (i)  If $\varphi : \bfa\to \bfb$ is a BP-map  then:
\pg (a)  $\kappa_{\bfb}= \hat\varphi (\kappa_\bfa) $ equivalently $\varphi$ is a BC-map of $(I,\kappa_\bfa)$ to $(J, \kappa_\bfb)$
\pg (b)  $\sigma_{\bfb}=\hat\varphi (\sigma_\bfa)$  equivalently $\varphi$ is a BC-map of $(I,\sigma_\bfa)$ to $(J, \sigma_\bfb)$
\p
\pg (ii) If  $\varphi : \bfa\to \bfb$ is an SP-map then:

 \pg (a) $\hat\varphi (\sigma_\bfa) \subset \sigma_{\bfb}$ equivalently $\varphi$ is an SC-map of $(I,\sigma_\bfa)$ to $(J, \sigma_\bfb)$
 \pg  (b) $\forall t\in \kappa_\bfb, \exists s\in \kappa_\bfa: \hat \varphi(s)\subset t$
 
\end{cor}

\begin{definota}\label{P-image}\rm
Let $\varphi : I \to J$ be a function. For any p-site $\bfa=\ppar{I, A, (A_i)}$ there exists one and only one p-site $\bfb$ with base $J$ and ground $A$ such that $\varphi$ is a BP-map of domain $\bfa$ and codomain $\bfb$,  namely the p-site $\ppar{J, A,  (B_j)}$ that satisfies relations (\ref{Bmap1}), equivalently the p-site with  parting function $\pi_\bfb :=  \hat \varphi \circ \pi_{\bfa}$. 
Such a p-site will be called  the \emph{P-image of $\bfa$ by $\varphi$} and is denoted $\varphi \bfa$. This is also the p-site  $\bfb$ with the finest profile $\bfb$  in $\mathrm {PROF}(J,A)$ that makes $\varphi$ an S-map from $\bfa$ to $\bfb$.  \hfill$\triangledown$
\end{definota}
In what follows we consider three arbitrary p-sites $\bfa, \bfb, \bfc$  with the same ground $A$, namely $\bfa\equiv \ppar{I,A, (A_i)}, \bfb\equiv \ppar{J,A, (B_j)}, \bfc\equiv \ppar{K,A, (C_k)}$.

\begin{pro} \label{Pcomposition} Let $\varphi : \bfa\to \bfb$ and $\psi: \bfb\to \bfc$ be BP-maps (resp. SP-maps)  then $\psi \circ \varphi$ is a BP-map  (resp. SP-map) with domain $\bfa$ and codomain $\bfc$.
\end{pro}
Proof:  For B-maps, we have for any $k\in K$, $C_k=\cup_{j\in \psi^{-1}(k)} B_j=$ $\cup_{j\in \psi^{-1}(k)}$ $(\cup_{i\in \psi^{-1}(j)}A_i)=$ $\cup_{i\in \varphi ^{-1}(\psi^{-1}(k))}A_i=\cup_{i\in (\psi\circ \varphi)^{-1}(k)}A_i$. For S-maps just replace sign $=$ by $\supset$ in the argument. \hfill $\Box$

\subsubsection{Categories of P-foundation  $\mathbf{bpol}_\bullet[\,\cdot\,, A]$, $\mathbf{spol}_\bullet[\,\cdot\,, A]$}

We introduce a category $\mathbf{bpol}_\bullet[\,\cdot\,, A]$,  (resp. ($\mathbf{spol}_\bullet[\,\cdot\,, A]$) whose objects  are  p-sites of ground $A$, and  whose morphisms of domain $\bfa$ and codomain $\bfb$  is the set of all B-maps (resp. the S-maps) denoted $\mathbf{bpol}_\bullet(\bfa, \bfb)$ (resp. $\mathbf{spol}_\bullet(\bfa, \bfb)$). The identity  on the the p-site $\bfa= \ppar{I, A, (A_i)}$ is defined as the B-map (resp. S-map) corresponding to the identity $id_I$. It is easily seen, in view of Proposition \ref{Pcomposition}, that the axioms of a category are satisfied. 
 \subsubsection{Partial Knit Functor $\kappa[\cdot, A]: \mathbf{bpol}_\bullet [\cdot, A]\to \mathbf{P}^2_=$}
If $\varphi$ is a BP- map of domain $\bfa$ and codomain $\bfb$ in  $\mathbf{bpol}_\bullet [\cdot, A]$,  then $\varphi$ is a BC-map  of domain $(I, \kappa_\bfa)$ and codomain $(J, \kappa_\bfb)$ in $\mathbf{P}^2_=$ (Proposition \ref{image2} (i) (a)). 
The association $\kappa :\bfa\mapsto (I, \kappa_\bfa)$  for objects and $\kappa : \varphi\mapsto \varphi$ for arrows satisfies the axioms relative to the image of composition laws.  It follows that $\kappa : \mathbf{bpol}_\bullet \to \mathbf{P}^2_=$ is a functor.

\begin{equation}
\begin{tikzcd}
 {\mathbf{P}^2_=}\\
{\mathbf{bpol}_\bullet[\cdot, A] }\arrow[u, "\kappa"]
\end{tikzcd}
\hskip12pt
\begin{tikzcd}[column sep=large]
{(I, \kappa_\bfa)}\arrow[r,  "{\varphi} " ]&{(J,\kappa_\bfb)}\\
{\bfa}\arrow[r, "{\varphi}" ]\arrow[u, maps to] &{\bfb}\arrow[u, maps to]
\end{tikzcd}
\end{equation}


\subsubsection{Partial Nerve Functor $\sigma[\cdot, A]: \mathbf{spol}_\bullet [\cdot, A]\to \mathbf{SIM}_\bullet$}\label{simcomplexes}

Now  we take as domain the category $\mathbf{spol}_\bullet$ instead of $\mathbf{bpol}_\bullet$.  $\sigma : \bfa \mapsto (I, \sigma_\bfa)$  for objects, and  $\sigma (\varphi) = \varphi$ for morphisms, is a functor of domain $\mathbf{spol}_\bullet$ and codomain $\mathbf{SIM}_\bullet$. 
Indeed, if $\varphi$ is an S-map from $\bfa$ to $\bfb$, then $\hat \varphi (\sigma_\bfa) \subset \sigma_\bfb$ (Proposition \ref{image2}) and the axioms for  composition laws are satisfied \ref{Pcomposition}.

\begin{equation}
\begin{tikzcd}
 {\mathbf{SIM}_\bullet}\\
{\mathbf{spol}_\bullet[\cdot, A] }\arrow[u, "\sigma"]
\end{tikzcd}
\hskip12pt
\begin{tikzcd}[column sep=large]
{(I, \sigma_\bfa)}\arrow[r,  "{\varphi} " ]&{(J,\sigma_\bfb)}\\
{\bfa}\arrow[r, "{\varphi}" ]\arrow[u, maps to] &{\bfb}\arrow[u, maps to]
\end{tikzcd}
\end{equation}

\subsection{P-foundation  with invariant base} \label{fixedsociety}
In his section all p-sites have the  same base denoted by $I$. We are going to consider the transformations that act on p-sites that keep the base invariant and that are induced by a \emph{shift} of the ground represented by a function. In fact, as in the preceding section, we define two classes of transformations that will lead to two categories of P-foundation with invariant base   $\mathbf{bpol}_\bullet[I,\cdot]$ and  $\mathbf{spol}_\bullet[I,\cdot]$.

\begin{defi}\label{map2} \rm  Let  $\bfa\equiv \ppar{I, A, (A_i)}$ and $\bfb \equiv\ppar{I,B, (B_i)}$  be two P-sites \emph{with the same base $I$}, and let  $f: A \to B$ be a function. We define two types of maps  with domain (or source)  $\bfa$ and codomain (or target) $\bfb$ .
\pg $f$ is a  \emph{BG-map}  if :
\begin{subequations}
\begin{equation} \label{Bmap}
A_i  =  f^{-1} (B_i)\,\,\, ( i\in I)
\end{equation}
$f$ is a  \emph{SG-map}  if :
\begin{equation} \label{Smap}
   f(A_i) \subset B_i\,\, \, (i\in I)
\end{equation}
\end{subequations}
\end{defi}

\subsubsection{BG-maps}

\begin{definota}\label{G-image}  \rm Let $f: A \to B$ a function. 
\pg For any p-site $\bfb \equiv\ppar{J, B, (B_j)}$, there exists one and only one p-site $\bfa \equiv\ppar{J, A, (A_j)}$ such that $f: \bfa\to \bfb$ is a BG-map. This is the site where $A_j= f^{-1}(B_j)$.  This p-site is called the \emph{inverse G-image} of $\bfb$ by $f$, and is denoted ${\bfb}f^{\!-\!1}$. It is characterized by its profile $\pi_\bfa = \pi_\bfb \circ f$.
\pg For any p-site $\bfa=\ppar{I, A, (A_i)}$, its \emph{direct G-image}  by  $f$ is the set $\bfb=\ppar{I,B, (B_i)}$ where $B_i=f(A_i)$ and it is denoted ${}^f\!\bfa$.  ${}^f\!\bfa$   has the finest cover in $\mathrm{cov}(I,B)$ that makes $f$ an S-map. \hfill$\triangledown$
\end{definota}

\begin{pro}\label{Bmap1}  Let $\bfa:= \ppar{I, A, (A_i)}$ and $\bfb :=\ppar{I,B, (B_i)}$ two P-sites. For any $f: A\to B$ the following are equivalent:
 \pg (i) $f$ is a BG-map with domain $\bfa$ and codomain $\bfb$
 \pg (ii) For any $s \in \calP(I)$   $f ^{-1}(\overset{\frown}B_s) =  \overset{\frown}A_s$
 \pg (iii) For any $s \in \calP(I)$   $f(\overset{\frown}A_s) \subset  \overset{\frown}B_s$
\pg  (iv) $\pi_\bfa^{}= \pi_\bfb ^{}\circ f$, in other terms, the following diagram is commutative:
\begin{equation}
\begin{tikzcd}[row sep=tiny]
{A}\arrow[dr,"\pi_\bfa"]\arrow[dd,"f" '] &{}\\
{}&\calP(I)\\
{B}\arrow[ur, "\pi_\bfb" ']&{}
\end{tikzcd}
\end{equation}
  
\end{pro}
Proof. $(i)\Leftrightarrow (ii)$  For a map $f: A \to B$ property  (\ref{Bmap}) implies $\overset{\frown}A_s = f^{-1} (\overset{\frown}B_s)$ for all $S \subset I$. Conversely  $A_s = \cup_{T\supset S} \overset{\frown}A_t =\cup_{T\supset S}  f^{-1}\overset{\frown}B_t = f^{-1}(\cup_{T\supset S} \overset{\frown}B_t)=f^{-1}(B_s)$, in particular for $S=\{i\}$ this is precisely (\ref{Bmap}).
\pg $(ii)\Rightarrow (iii)$ is straightforward.

\pg   $(iii)\Rightarrow (iv)$ : 
 For any  $a \in A$, put $S:= \pi_\bfa (a)$,  so that  $a \in \overset{\frown}A_s$ and  $f(a) \in  f(\overset{\frown}A_s) \subset \overset{\frown}B_s$. Thus $\pi_\bfb^{} (f(a)) =S=  \pi_\bfa (a)$.  
 \pg $(iv)\Rightarrow (ii)$ :  For any $S\subset I$, $\overset{\frown}A_s= \pi_\bfa ^{-1} (S) =  f^{-1} (\pi_\bfb ^{-1} (S)) = f^{-1}(\overset{\frown}B_s)$.

\hfill $\Box$

\begin{fact} \label{compositionBG} Let $ f : \bfa\to \bfb$ and $ g:  \bfb\to \bfc$ be BG-maps   then $g \circ f$ is a BG-map  with domain $\bfa$ and codomain $\bfc$.
\end{fact}
Proof:  Straightforward \hfill $\Box$

 \begin{pro} Let $f: A \to B$.  Then $f$ is an isotopy  of $\bfa$ to $\bfb$ if and only if $f$ is a bijective BG-map. 
 \end{pro}

In fact isotopy is too strong  a requirement to express sameness of p-sites in our framework. Below we shall define the notion of BG-isomorphism which turns out to be  adapted to our study.

\begin{pro} \label{Bexistence}There exists  a BG-map $f: \bfa \to \bfb$  if and only if $\kappa_\bfa \subset \kappa_\bfb$.  
\end{pro}
 Proof. Assume that $f$ is a B-map. In view of Lemma \ref{Bmap1} (iii), for any $S\subset I$,  $\overset{\frown}A_s \neq \emptyset$ implies $\overset{\frown}B_s \neq \emptyset$ so that  $ \kappa_\bfa \subset \kappa_\bfb$. 
Conversely  first select  for any $S \in \kappa_\bfb$ an element $b_s \in B_s$.  If  $\kappa_\bfa \subset \kappa_\bfb$,  one can define a map on $A$ by positing for any $a \in \overset{\frown}A_s$  $f(a): =b_s$.   In view of Lemma \ref{Bmap1} (ii), $f$ is clearly a BG-map. \hfill $\Box$

\begin{defi} \label{Bisomorphism}\rm  Two p-sites $\bfa$ and $\bfb$ with the same base $I$, are said to be \emph{BG-isomorphic} if there is a BG-map $f: \bfa\to \bfb$ and a BG-map  $g: \bfb\to \bfa$. \hfill$\triangledown$
\end{defi}

\g It follows from the above  that $\bfa$ and $\bfb$ are  BG-isomorphic if and only if 
 $\kappa_\bfa = \kappa_\bfb$ 
 
 \begin{fact} \label{lem3} Let $\bfa\equiv\ppar{I, A, (A_i)}$ and $\bfb \equiv \ppar{I, B, (B_i)}$   be two p-sites,  $f : \bfa \to \bfb$ a \emph{surjective}  BG-map.  There exists a right inverse  $g: B \to A $ to  the function $f$ ( $f\circ g=id_B$). Any right inverse  to $f$ induces a BG-map. In particular $\bfa$ and $\bfb$ are BG-isomorphic. 
\end{fact}
\emph{Proof}. If $f\circ g =id_B$, then $B_i = g^{ -1} ( f^{-1} (B_i))=g^{-1}(A_i)$. \hfill$\Box$

\subsubsection{SG-maps}
\begin{pro}\label{Smap1}  Let $\bfa \equiv \ppar{I, A, (A_i)}$ and $\bfb\equiv \ppar{I,B, (B_i)}$ two p-sites with the same base $I$.
  For any  function $f: A\to B$ the following are equivalent:
 \pg (i) $f$ is an SG-map with domain $\bfa$ and codomain $\bfb$
 \pg (ii) For any  $s\in \calP_0(I))$,  $f( A_s) \subset B_s$
 \pg (iii)  For any  $s\in \calP_0(I))$,  $A_s \subset f^{-1} ( B_s) $
\pg  (iv) $\pi_\bfa \subset \pi_\bfb \circ f $
\end{pro}
\pg Proof.  \pg   $(i)\Leftrightarrow (ii) \Leftrightarrow (iii)$   is straightforward.
\pg  \pg  $(iii)\Rightarrow (iv)$ : For any  $a \in A$, put $s:= \pi_\bfa (a)$,  so that  $a \in \overset{\frown}A_s \subset A_s$. Thus  $f(a) \in B_s  = \sqcup\{\overset{\frown}B_t: t\supset s\}$. It follows that  $f(a) \in \overset{\frown}B_t$ for some $t \supset s$ and therefore $\pi_\bfb(f(a)=t  \supset s = \pi_\bfa (a). $
\pg $(iv)\Rightarrow (iii)$. Let $a \in A_s$. Since $A_s = \sqcup\{\overset{\frown}A_t: t\supset s\}$ it follows that  $\pi^{}_\bfa (a) \supset s$,   and  consequently $  \pi^{}_\bfb (f(a))\supset \pi^{}_\bfa (a) \supset s $. We conclude that  $f(a) \in B_s$. \hfill $\Box$


\begin{fact} \label{compositionBG} Let $ f : \bfa\to \bfb$ and $ g:  \bfb\to \bfc$ be SG-maps   then $g \circ f$ is an SG-map  with domain $\bfa$ and codomain $\bfc$.
\end{fact}
Proof:  Straightforward \hfill $\Box$

 \begin{pro} Let $f: A \to B$.  Then $f$ is an isotopy  of $\bfa$ to $\bfb$ if and only if $f$ is a bijective and  both $f$ and its inverse are  SG- maps.
 \end{pro}

\pg Below we shall define  the notion SG-isomorphism which is adapted for the study of political structures (simplicial complexes).

\begin{pro} \label{Sexistence}There exists  a  SG-map $f: \bfa \to \bfb$  if and only if $\sigma_\bfa\subset \sigma_\bfb$.  
\end{pro}
 \emph{Proof}. 
Assume that $f$ is an S-map. In view of Lemma \ref{Smap1}, for any $S\subset I$,  $ A_s \neq \emptyset$ implies $B_s \neq \emptyset$, therefore $\sigma_\bfa\subset \sigma_\bfb$.  Conversely assume the inclusion $\sigma_\bfa\subset \sigma_\bfb$.  We start by choosing  $b_s\in B_s$ for all $S\in \sigma_\bfb$, and an arbitrary element $b_\emptyset \in B$.  Since    ${\kappa_0} _\bfa \subset  \sigma_\bfa$ we have  ${\kappa_0} _\bfa \subset  \sigma_\bfb$. Therefore we can define  for any $S\in {\kappa_0}_\bfa$,  $f$ on $\overset{\hskip2pt\frown}A_s$  by $\varphi (a)= b_S$, and in case $\emptyset \in \kappa_\bfa$ by $ f(a) = b_\emptyset$ on  $\bar A_\emptyset$. We see that,  for $s\in \calP_0(I)$,  $f(A_s) = f (\cup_{t\supset s} \overset{\hskip2pt\frown}A_t)= \cup_{t\supset s}f(\overset{\hskip2pt \frown}A_t)=\cup_{t\supset s} \{B_t\}\subset B_s$.     In view of Lemma \ref{Smap1} $f$ is an S-map.  \hfill $\Box$

\begin{defi} \label{Sisomorphism}\rm  Two p-sites $\bfa$ and $\bfb$ with the same base $I$, are said to be \emph{SG-isomorphic} if there is an SG-map $f: \bfa\to \bfb$ and an SG-map  $g: \bfb\to \bfa$. \hfill$\triangledown$
\end{defi}
 
 \g It follows from the above  that $\bfa$ and $\bfb$ are  S-isomorphic if and only if 
 $\sigma_\bfa = \sigma_\bfb$ 
 \begin{example} (Gallopolis continued)\label{gallopolis2}
Following the 2024 parliamentary elections in Gallopolis (Example \ref{gallopolis}) the new distribution of party strength has rendered the previously governing coalition -- comprising CONS and  LiBER  --  insufficient to secure a majority.  In response to this deadlock, a bargaining process was initiated to explore the formation of a new coalition.
 We assume that this process  unfolded around three hypothetical scenarios, each defined by the omission of one of the three foundational dimensions that structure the political space. The aim  is to assess whether viable and winning coalitions can emerge under reduced dimensional ground.

\pg \emph {Scenario 1 : Disregarding the Openness dimension}. This scenario preserves the political configuration that prevailed prior to the elections (see Graphic \ref{figgallopolis1}). However  no  viable   and winning coalition exists in the newly elected parliament. Consequently this scenario replicates the current impasse and fails to offer a viable path for governance.

\pg \emph {Scenario 2 :  Disregarding the Economic dimension}. In this scenario  (see graphic \ref {figgallopolis2})  the  configuration allows the formation of a broad  coalition  comprising  LIBER, SOCD, LEFT. This coalition is viable since the state  $(l, \gamma)$ lies in the intersection  of their projected aspiration sets.  Importantly  this coalition constitutes a winning majority under the new parliament composition, suggesting that the omission of economic cleavages could enable effective coalition-building on the base of shared societal and openness values.  

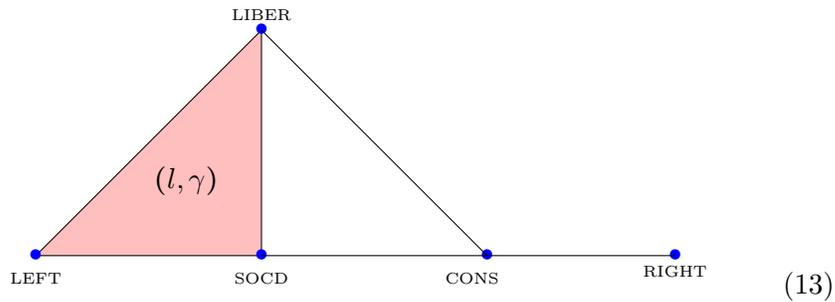
\begin{center}
\begin{equation}\label{figgallopolis2}
\begin{tikzpicture}[scale=1]
\draw[fill=pink] (0,0)--(3,0)--(3,3)--cycle;
\draw[blue] (0,0) node {$\bullet$};
 \draw (0,-0.3)node  {\tiny LEFT};
\draw (2,1)node {$(l,\gamma)$};
\draw[blue] (3,0) node{$\bullet$}; \draw  (3, -0.3) node {\tiny SOCD};
\draw[blue] (6,0) node{$\bullet$}; \draw (5.8,-0.3) node {\tiny CONS};
\draw[blue] (3, 3) node{$\bullet$} node[above, black] {\tiny LIBER};
\draw[blue] (8.5,0) node{$\bullet$} node[below, black] {\tiny RIGHT};
\draw(3,0)--(6,0);
\draw(6,0)--(8.5,0);
\draw(3,3)--(6,0);

\end{tikzpicture}
\end{equation}
\end{center}
  
\pg \emph{Scenario 3 : Disregarding the Societal dimension}. In this scenario (see  Graphic \ref{figgallopolis3})   the coalition comprising LIBER, CONS, RIGHT  is viable since the aspirations of the three parties  agree on state $(1, \alpha)$, and it is  winning  in the new parliament.

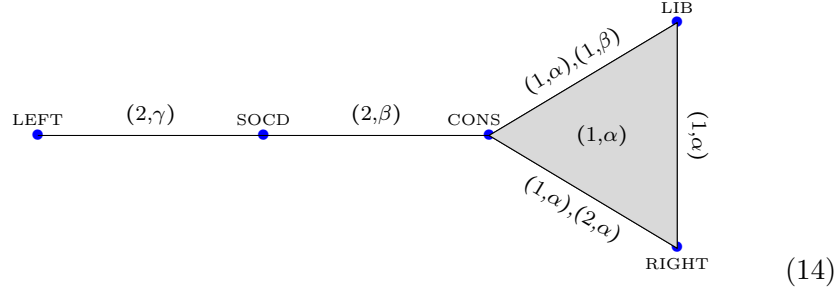
\begin{center}
\begin{equation}\label{figgallopolis3}
\begin{tikzpicture}[scale=1]
\draw[blue] (0,0) node{$\bullet$} node[above, black] {\tiny LEFT};
\draw[blue] (3,0) node{$\bullet$} node[above, black] {\tiny SOCD};
\draw[blue] (6,0) node{$\bullet$}; \draw(5.8,0) node[above, black] {\tiny CONS};
\draw[blue] (8.5,1.5) node{$\bullet$} node[above, black] {\tiny LIB};
\draw[blue] (8.5,-1.5) node{$\bullet$} node[below, black] {\tiny RIGHT};
\draw(0,0)--(3,0)node[midway,above]{$\scriptstyle (2, \gamma)$}--(6,0)node[midway,above]{$\scriptstyle (2, \beta)$};
\draw(6,0)--(8.5,1.5)node[midway,above,sloped] {$\scriptstyle (1,\alpha), (1,\beta)$};
\draw(6,0)--(8.5,-1.5)node[midway,below,sloped] {$\scriptstyle (1,\alpha), (2, \alpha)$};
\draw(8.5,1.5)--(8.5,-1.5)node[midway, above,sloped] {$\scriptstyle (1,\alpha)$};
\draw[fill=gray!30](6,0)--(8.5,1.5)--(8.5,-1.5)-- cycle;
\draw(7.5, 0) node{$\scriptstyle  (1,\alpha)$};
\end{tikzpicture}
\end{equation}
\end{center}

Note that under the Gallopolis Constitution the President has  the last word in the decision process. Scenario 2 and Scenario 3 were available, since both have winning coalitions. After a bargaining process that lasted three months, the President chose scenario 3, a coalition on the Right.

  Disregarding a dimension $x$ of the Ground space amounts,  in mathematical terms, to a shift of the ground of p-site $\bfa= \ppar{I, A, (A_i)} $.  Such a transformation is induced by a the projection  say $p^x$ of $A$ onto $B$ parallel to the disregarded dimension.  The resulting p-site is the G-image of $\bfa$ by $p^x$ denoted ${}^{p^x}\!\bfa \equiv \ppar{I, B, (B_i)}$ where $B=p^x(A)$ and $B_i= p^x(A_i)$ (see definition and notation \ref{G-image}). 
   \hfill $\bigcirc$

\end{example}

\subsubsection{Category  of P-foundation $\mathbf{bpol}_\bullet[I,\cdot]$} \label{bpol}Objects of this category are  p-sites and the morphisms of domain $\bfa$ and codomain $\bfb$, denoted $\mathbf{bpol}_\bullet[I, \cdot](\bfa, \bfb)$.The identity on $A$  induces an identity on $\bfa$. In view of Proposition \ref{compositionBG} the law of composition of morphisms is satisfied. 

\subsubsection{Functors $\kappa[I, \cdot]$ (Knit) and $\bca[I]$ (Canon) } 
There is a functor $\kappa\equiv \kappa[I, \cdot]: \mathbf{bpol}_\bullet[I,\,\cdot\,] \to \mathrm{P} ^2_\subset[I] $ defined by the association  $\bfa\to \kappa_\bfa$ and for  a BG-map $f: \bfa\to \bfb$,  $\kappa (f)$ is the inclusion $\kappa_\bfa \hookrightarrow\kappa_\bfb$. It is easy to verify the functor axioms   by relying on Proposition \ref{compositionBG} and Proposition \ref{Bexistence}.

\begin{equation}
\begin{tikzcd}
 {\mathbf{P}^2_\subset[I]}\arrow[d, "\bca",  xshift =1ex]\\
{\mathbf{bpol}_\bullet[I, \cdot] }\arrow[u, "\kappa", xshift=-1ex]
\end{tikzcd}
\hskip12pt
\begin{tikzcd}[]
{\kappa_\bfa}\arrow[r, hook, " i^{\kappa_\bfa}_{\kappa_\bfb}" ]&{\kappa_\bfb}\\
{\bfa}\arrow[r,maps to," f " ']\arrow[u," " '] &{\bfb}\arrow[u, maps to] 
\end{tikzcd}
\hskip12pt
\begin{tikzcd}[]
{\calX}\arrow[r,hook,"  i^{\calX}_{\calY} "]\arrow[d, maps to] &{\calY}\arrow[d, maps to] \\
{\bca(\calX)}\arrow[r,maps to, " i^{\calX}_{\calY}" ']&{\bca(\calY)}
\end{tikzcd}
\end{equation}

\p
On the opposite direction we define a functor $\bca\equiv\bca[I]: \mathbf{P} ^2_\subset(I) \to  \mathbf{bpol}_\bullet[I,\,\cdot\,]$ by  associating  to any object $X \in \calP^2(I)$  the canonical p-site associated to $\calX$, namely  $\bca(\calX) := \bca (I,\calX)$ (Definition \ref{canonicaldefi}), and to any morphism $i^\calX_\calY$ we put   $\bca (i^\calX_\calY) := {i}^\calX_\calY$. Indeed   if $i^\calX_\calY$ is a morphism of $\mathrm{ P_\subset ^2(I)}$ one has  $\pi_{\bca(\calX)}= id_\calX = id_\calY \circ i^\calX_\calY = \pi_{\bca(\calY)} \circ i^\calX_\calY$.  It follows, in view of Proposition \ref{BGcharacter}, that   the inclusion  $i^\calX_\calY$ is  a  BG-map between $\bca(\calX)$ and $\bca (\calY)$.  It is clear that $\bca$ is a functor.
\subsubsection{Structure relations between $\mathbf{bpol}_\bullet[I, \cdot]$ and $\mathbf{P}^2_\subset(I)$, a sketch} 
\pg We search for  an ``isomorphism'' between $\mathbf{bpol}_\bullet[I, \cdot]$ and $\mathbf{P}^2_\subset(I)$ via an ``isomorphism'' between the functor   $\bca \circ \kappa :\mathbf{bpol}_\bullet[I, \cdot] \to  \mathbf{bpol}_\bullet[I, \cdot]$  and the functor  identity $\mathbf{1}_{\mathbf{bpol}_\bullet[I, \cdot]}$. It is clear  that  $\bca \circ \kappa$ is not equal  to $\mathbf{1}_{\mathbf{bpol}_\bullet[I, \cdot]}$, that is the functors $\bca$ and $\kappa$  are not  strictly speaking inverse of each other.  In terms of category theory one rather seeks to show a weaker type of isomorphism, namely  that  nevertheless  $\bca \circ \kappa$ is \emph{naturally isomorphic} to  $ \mathbf{1}_{\mathbf{bpol}_\bullet[I, \cdot]}$. 
The following sketch of the method in the present case  provides  a template for the general figure postponed to Section 5 
\pg \emph{Step 1. Natural transformation.}  First  we show   that there is a natural transformation that takes $\mathrm{1}_{\mathbf{bpol}_\bullet[I, \cdot]}$ to $\bca \circ \kappa$. 
The ground of  $ \bca\circ \kappa (\bfa)\equiv \bca(I, \kappa_\bfa)$  is  the set $\kappa_\bfa$. 
\pg
Since  $\overline{\pi}_\bfa \in \mathbf{bpol}_\bullet[I, \cdot](\bfa, \bca\circ \kappa (\bfa))$ ( Lemma \ref{pi-Bmorphism})  we are  led to consider the transformation between the functors $\mathrm{1}_{\mathbf{bpol}_\bullet[I, \cdot]}$ to $\bca \circ \kappa$, provided by 
 the family of morphisms  :  $\overline{\pi}= (\overline{\pi}_\bfa)$. This provides  the fundamental structure result: 
\begin{pro}\label{transformationB} $\overline{\pi}$  is a natural transformation from $\mathrm{1}_{\mathbf{bpol}_\bullet[I, \cdot] }$  to $ \bca\circ \kappa$ also written in symbols  $ \mathrm{1}_{\mathbf{bpol}_\bullet[I, \cdot] }\overset{\pi} \Longrightarrow \bca\circ \kappa$ 
\end{pro}

\pg \emph{Step 2.  Quotient on the morphisms.} Let  $\mathbf{bpol}_{\bullet/b}[I,\,\cdot\,]$ be the category that has the same objects as  $\mathbf{bpol}_\bullet[I,\,\cdot\,]$ and where there is a morphism between $\bfa$ and $\bfb$ if  the family  of all B-maps between $\bfa$ and $\bfb$ is nonempty.  The identity of $\bfa$  is thus the family of all  functions  $f: A\to A$ such that $\pi_\bfa \circ f = \pi_\bfa$.  An isomorphism between $\bfa$ and $\bfb$  in this category exists if and only if there is a B-map between $\bfa$ and $\bfb$ and a B-map between $\bfb$ and $\bfa$.   We prove the first structure theorem: 
 
\begin{thm}\label{thm1B}
$\begin{tikzcd}
{\mathbf{bpol}_{\bullet/b}[I, \cdot]}\arrow[r, "{\kappa} ", yshift =0.7ex]&{\mathbf{P}^2_\subset [I]}\arrow[l, " \bca"  , yshift=-0.7ex]
\end{tikzcd}$
is an equivalence of categories
\end{thm}

\pg \emph{Step 3. Quotient on the objects. } The class of p-sites is partially preordered  by the relation : $\bfa \to_B\bfb$ if and only if  $\mathbf{bpol}_\bullet[I, \cdot][\bfa, \bfb]$ is non-empty. The equivalence relation $\bfa \sim_B \bfb$ if and only if $\bfa$ and $\bfb$ are B-isomorphic (Definition \ref{Bisomorphism}) induces a poset ${\mathbf{BPOL}_\bullet[I, \cdot]}$. We prove the second structure theorem :
 
 \begin{thm} \label{Bposet1}  The poset $\mathbf{BPOL}_\bullet[I, \cdot]$ is isomorphic  to  $\mathbf{P}^2_\subset[I]$. 
\end{thm}

\section{P-foundation with variable Base and Ground}
In subsection \ref{fixedground} the ground is invariant and the base is variable, in subsection \ref{fixedsociety} the base is invariant and the ground is variable. 
In this section  we extend the definitions and the results  to the  general model in which both the base and the ground are variable. The proofs follow the same pattern as in the above section and therefore repetitions  in the proofs can be avoided, but for clarity some details  are maintained. 
\subsection{B-maps, S-maps}
\begin{defi} \rm  Let $\bfa:= \ppar{I, A, (A_i)}$ and $\bfb :=\ppar{J,B, (B_j)}$  be two p-sites $\varphi: I\to J$,   $f: A \to B$ be two functions. We define two classes of maps  with domain (or source)  $\bfa$ and codomain (or target) $\bfb$ .
\pg $(\varphi,f)$ is a \emph{B-map}  if :
\begin{equation}\label{Bmap2}
 \underset{i\in \varphi ^{-1}(j)}\cup A_i  =  f^{-1} (B_j)\,\,\, ( j\in J)
\end{equation}
\pg $(\varphi,f)$ is a  \emph{S-map} if :
\begin{subequations}
\begin{equation}\label{Nmap2}
 \underset{i\in \varphi^{-1}(j)}\cup A_i \subset  f^{-1} (B_j)\,\, \, (j\in J)
 \end{equation}
 or equivalently : 
 \begin{equation}
 f(i) =j \Rightarrow f(A_i) \subset B_j
 \end{equation}
\end{subequations}
\hfill$\triangledown$
\end{defi}


Let $\bfa:= \ppar{I, A, (A_i)}$ and $\bfb :=\ppar{J,B, (B_j)}$  be two P-sites $\varphi: I\to J$,   $f: A \to B$ be two functions.
We recall the P-image  of $\bfa$ by $\varphi$  denoted $\varphi\bfa := \ppar{J, A, \hat\varphi\circ\pi_\bfa}$ (Notation and definition \ref{P-image})  and  the inverse G-image of $\bfb$  by $f$  denoted $\bfb f^{\!-1}:= \ppar{J, A,( f^{\!-1}(B_j)} \equiv \ppar{J, A, \pi_\bfb\circ f} $ (Notation and definition \ref{G-image}).  We have the following sequence of functions :  

\begin{equation}
(I,A)\, \overset{\varphi, id_A} \longrightarrow\, (J,A)\,\overset{id_J, id_A}\longrightarrow\, (J,A) \overset{id_J, f} \longrightarrow\, (J,B) 
\end{equation}
and the corresponding  sequence of  maps and p-sites:
\begin{equation}
\bfa\, \overset{\varphi, id_A} \longrightarrow\, \varphi \bfa\,\overset{id_J, id_A}\longrightarrow\, \bfb f^{\!-\!1} \overset{id_J, f} \longrightarrow\,\bfb 
\end{equation}
Note that the first map $\bfa\, \overset{\varphi, id_A} \longrightarrow\, \varphi\bfa$ and the last one $\bfb f^{\!-1} \overset{id_J, f} \longrightarrow\,\bfb$  are  B-maps. The nature of the sequence is therefore determined by the central map $\varphi \bfa\,\overset{id_J, id_A}\longrightarrow\, \bfb f^{\!{-\!1}}$. The following provides equivalent properties in each case:

 \begin{pro}\label{BScharacter} Let  $\varphi: I \to J$ and $f: A \to B$  be two functions. Then  the following are equivalent: 
 \pg (i) $(\varphi,f) : \bfa \to \bfb$ is a B-map (resp. S-map) 
 \pg (ii) $( id_J, f):  \varphi\bfa\to \bfb$ is a B-map (resp. S-map)
 \pg (iii) $(\varphi, id_A) : \bfa\to {\bfb}f^{\!-\!1}$ is a B-map (resp. S-map)
\pg (iv) $(id_J, id_A): \varphi\bfa\to \bfb f^{\!-\!1}$ is a B-map (resp. S-map)
 \pg (v) $\hat \varphi \circ \pi_\bfa= \pi_\bfb \circ f$  (resp. $ \hat \varphi \circ \pi_\bfa\subset\pi_\bfb \circ f$)
\end{pro}
\emph{Proof}.  The B-case is a consequence of Proposition \ref{Bmap1} and Proposition \ref{image2}, the S-case is a consequence of Proposition \ref{Smap1} and Proposition \ref{image2} \hfill $\Box$
\pp
In order to define  formally the categories  of the next subsections we need the following:
\begin{pro} \label{category1}Given the sequence of  B-maps (resp. S-maps) $\bfa \overset{\varphi,f}\rightarrow\bfb\overset{\psi,g}\rightarrow\bfc$ 
 their composite defined by  $(\psi,g)\circ (\varphi,f):= (\psi\circ\varphi, g\circ f)$ is  a B-map (resp. S-map) $\bfa\to \bfc$.
\end{pro}
\emph{Proof}.  The B- case. $\pi_\bfc\circ(g\circ f)=(\pi_\bfc\circ g)\circ f=(\hat\psi\circ\pi_\bfb) \circ f=\hat\psi\circ(\pi_\bfb \circ f) =\hat\psi \circ (\hat \varphi \circ \pi_\bfa) = (\hat\psi \circ\hat \varphi) \circ \pi_\bfa= \widehat {(\psi\circ\varphi)} \circ \pi_\bfa$.
\pg The S-case  is similar. Just replace the second and forth equality sign $=$  by $\supset$. \hfill$\Box$  
\begin{fact}\label{category2} Given the sequence of  B-maps (resp. S-maps) $\bfa \overset{\varphi,f}\rightarrow\bfb \overset{\psi,g}\rightarrow\bfc \overset{\rho,h}\rightarrow \bfd $ one has:
$\big((\rho,h)\circ \psi,g)\big)\circ (\varphi,f) = (\rho,h)\circ \big(\psi,g)\circ (\varphi,f)\big)$  (Associativity)
\end{fact}
 \begin{pro}\label{BSexistence2}  Let $\varphi : I\to J$ be a function and let $\bfa$ and $\bfb$ be p-sites. There exists $f : A\to B$ such that $(\varphi, f )$  is a B-morphism  (resp. S-morphism) with domain $\bfa$ and codomain $\bfb$ if and only if $\hat\varphi(\kappa_\bfa)\subset \kappa_{\bfb}$  (resp. $\hat \varphi(\sigma_{\bfa}) \subset \sigma_{\bfb}$)
  \end{pro}
\emph{Proof}. We prove the B-case. Recall first that $\pi_{\varphi \bfa} = \hat\varphi(\kappa_\bfa)$ (Definition and notation \ref{P-image}.  If   $\pi_{\varphi \bfa}\subset \kappa_{\bfb}$ then in view of Proposition \ref{Bexistence} there exists a function $f: A\to B$ such that $(id_J, f)$ is a B-map from $\varphi\bfa$ to $\bfb$. Now in view of  Proposition \ref{BScharacter} this implies that $(\varphi, f)$ is a B-map from $\bfa$ to $\bfb$. The proof of the B-case is similar, with reference to Proposition \ref{Sexistence}  instead of  Proposition \ref{Bexistence}. \hfill $\Box$
\pp
The following extends the definition of B-isomorphism to  arbitrary p-sites:
\begin{defi} \label{BSisomorphism}\rm  Two p-sites $\bfa\equiv \ppar{I, A, (A_i)}$ and $\bfb\equiv \ppar{J,B, (B_j)}$ are said to be \emph{B-isomorphic} (resp.  \emph{S-isomorphic} ) if there are B-maps (resp. S-maps) $(\varphi, f) :\bfa\to \bfb$  and $(\psi, g) : \bfb\to \bfa$  such that $\psi\circ\varphi =id_I$ and $\varphi\circ \psi =id_J$. \hfill$\triangledown$
\end{defi}

The following will be useful for the definition of functors of the next subsections.
\begin{pro} Let $\calX\in \calP^2(I)$, $\calY\in \calP^2(J)$,  $\varphi : I\to J$, $f: \calX\to \calY$  two functions.   $(\varphi, f)$ is a B-map $\bca(I, \calX)$ $\to$ $\bca(J, \calY)$ if and only if $f(s)=\hat\varphi (s)$ for all $s\in \calX$. In particular the condition implies $\hat \varphi (\calX)\subset \calY$, or said otherwise, that $\varphi$ is an SC-map $(I,\calX) \to (J, \calY)$
\end{pro}
\emph{Proof}. Put $\bfa:= \bca(I, \calX)$ and $\bfb :=\bca(J, \calY)$. $(\varphi, f)$ is a B-map  if and only if $\hat \varphi\circ \pi_\bfa= \pi_\bfb\circ f$. Since $\pi_\bfa =i^\calX$ and $\pi_\bfb=i^\calY$, the condition is equivalent to $\hat \varphi (s) =f(s)$ for all $s\in \calX$. This implies in particular $\hat \varphi (s)\in \calY$ for all $s\in \calX$. \hfill $\Box$

\subsection{ Categories of P-foundation versus P-configuration}
This subsection is devoted to the study of the structural relationships between the P-foundation categories and the P-configuration categories. 
Our model reveals  two theories of P-foundation, a B-theory where the morphisms are built on B-maps and an S-theory where the morphisms are built on S-maps. In each one we give the formal definition of the category and we compare it via appropriate functors to the P-configuration category.
\subsubsection{Category $\mathbf{bpol}_\bullet$  }\label{pol-B}
Let $\mathbf{bpol}_\bullet [\cdot, \cdot ]$  be the category whose objects are p-sites $\bfa, \bfb, \ldots$  and whose morphisms between $\bfa$ and $\bfb$, denoted $\mathbf{bpol}_\bullet(\bfa,\bfb)$, are B-maps $(\varphi, f)$ with domain $\bfa$ and codomain $\bfb$. The identity on $\bfa$ is defined as the ordered pair $(id_I, id_A)$. In view of Proposition \ref{category1} and Fact \ref{category2} that the category axioms are indeed verified. In what follows, we put   $\mathbf{bpol}_\bullet$ in correspondence with the category  $\mathbf{P}^2_\bullet$  introduced in subsection \ref{complexes}.


\paragraph{The Knit functor $\scrK$} --
If $(\varphi,f)$ is a B-map of domain $\bfa$ and codomain $\bfb$, it follows from \ref{BScharacter}(v) that $\hat\varphi (\kappa_\bfa) \subset \kappa_\bfb$, and therefore $\varphi$ is a morphism of domain $(I, \kappa_\bfa)$ and codomain $(J,\kappa_\bfb)$. 
The association $\scrK:\bfa\mapsto (I, \kappa_\bfa)$  for objects and $\scrK : (\varphi, f) \mapsto \varphi$ for arrows is  well defined  and is represented by the diagram  
\begin{equation}
\begin{tikzcd}
 {\mathbf{P}^2_\bullet}\\
{\mathbf{bpol}_\bullet }\arrow[u, "\scrK"]
\end{tikzcd}
\hskip12pt
\begin{tikzcd}[column sep=large]
{(I, \kappa_\bfa)}\arrow[r, hook, "{\varphi} " ]&{(J,\kappa_\bfb)}\\
{\bfa}\arrow[r, "{\varphi, f}" ]\arrow[u, maps to] &{\bfb}\arrow[u, maps to]
\end{tikzcd}
\end{equation}

\pg and has  the following properties:
\pg (i) for every $\bfa = \ppar{I,A, (A_i)}$, $\kappa(id_I, id_A)=id_I$
\pg (ii)  $\kappa(\psi\circ\varphi, g\circ f) = \psi \circ \varphi=\kappa(\psi, g) \circ \kappa(\varphi, f)$ 

\pg It follows that $\scrK : \mathbf{bpol}_\bullet \to \mathbf{P}^2_\bullet$ is a functor. 
\p
\paragraph{The Functor Canon $ \scrC$} -- In the opposite direction the association $\scrCa: (I,\calX) \mapsto \bca(I,\calX)$  for objects and $\scrCa : \varphi \mapsto (\varphi, {\hat\varphi}^X_Y )$  for morphisms is represented by the diagram

\begin{equation}
\begin{tikzcd}
 {\mathbf{P}^2_\bullet}\arrow[d, "\scrCa"]\\
{\mathbf{bpol}_\bullet }
\end{tikzcd}
\hskip12pt
\begin{tikzcd}[column sep=large]
{(I,\calX)} \arrow[d, maps to]\arrow[r,hook,   "{\varphi} " ]&{(J,\calY)}\arrow[d, maps to]\\
{\bca(I,\calX)}\arrow[r,  "{\varphi,{\hat\varphi}^\calX_\calY} "] &{ \bca(J,\calY)} 
\end{tikzcd}
\end{equation}
Precisely if $\varphi: (I,\calX)\to (J,\calY)$ is a morphism of p-formations,  one has  $\pi_{\bca(J,\calY)}\circ \hat\varphi^X_Y = \hat \varphi \circ  \pi_{\bca(I,\calX)}$.  It follows from Proposition \ref{BScharacter}(v), that  the couple $(\varphi, \hat\varphi^X_Y)$ is  a  B-map between $\bca(I,\calX)$ and $\bca (J,\calY)$.
It is easy to verify that $\scrC : \mathrm{P}^2_0 \to \mathbf{bpol}_\bullet$  is a functor.
\paragraph{The natural transformation $\overline{\pi}$} -- 
\begin{lem}\label{pi-Bmorphism} For any p-site $\bfa$,  the parting  function of $\bfa$  considered as a map $A \to \kappa_\bfa$, denoted $\overline{\pi}_\bfa$,  induces  a BG-map  $\bfa \to \bca(I, \kappa_\bfa)$ 
\end{lem}
\emph{Proof}.  Put $\bfb:= \bca(I, \kappa_\bfa)$,  then  $ \pi_\bfb =i^{\kappa_\bfa}_{\calP(I)}$, Then one has  $\pi_\bfa= \pi_\bfb\circ \overline{\pi}_\bfa$ therefore in view of Proposition \ref{Bmap1}(iv)), $\overline{\pi}_\bfa$ is a BG-map. \hfill $\Box$

\pg
In view of  Lemma \ref{pi-Bmorphism},  $\overline{\pi}_\bfa$ is  a morphism. We are going to consider the family of morphisms  in $ \mathbf{bpol}_\bullet[I, \cdot]$:  $\overline{\pi}= (\overline{\pi}_\bfa)$  as a map of  domain the functor  $\mathrm{1}_{\mathbf{bpol}_\bullet}$ and codomain the functor $ \scrCa\circ \scrK$.

\begin{lem} \label{transformationB2} $\overline{\pi}$  is a natural transformation from $\mathrm{1}_{\mathbf{bpol}_\bullet}$  to $ \scrCa\circ \scrK$ also written in symbols  $\mathrm{1}_{\mathbf{bpol}_\bullet} \overset{\overline \pi}\Longrightarrow \scrCa\circ \scrK$. 
\end{lem}
\emph{Proof.} For every pair $\bfa$, $\bfb$ of p-sites and every morphism $(\varphi, f):\bfa\to \bfb$, the following diagram is commutative:

\begin{equation}
\begin{tikzcd}[column sep=large]
{\bfa}\arrow[d, "{\overline\pi}_\bfa" ']\arrow[r, "{\varphi, f}" ] &{\bfb}\arrow[d, "{\overline\pi}_\bfb "] \\
{\scrCa\circ \scrK(\bfa)}\arrow[r,  " {\scrCa\circ \scrK (\varphi, f)}" ']&{\scrCa\circ \scrK(\bfb)}
\end{tikzcd}
\end{equation}
One has only to remark that $\scrCa\circ \scrK(\bfa) = \bca(I, \kappa_\bfa)$, $\scrCa\circ \scrK(\bfb)= \bca(I, \kappa_\bfb)$, and the lower map $\scrCa\circ \scrK (\varphi, f)$ is the B-map $(\varphi, \hat\varphi^{\kappa_\bfa}_{\kappa_\bfb}): \bca(I, \kappa_\bfa)\to \bca(J, \kappa_\bfb). $\hfill $\Box$

\paragraph{$\mathbf{bpol}_{\bullet/b}$ and Category equivalence to $\mathbf{P}^2_\bullet$} \label{POL-B}
 -- The category $\mathbf{bpol}_{\bullet/b}$ is defined as follows:
\pg (i) The objects are p-sites,
\pg (ii) The set of morphisms of domain $\bfa$ and codomain $\bfb$ is : 
\pg  $\mathbf{bpol}_{\bullet/b}(\bfa, \bfb) : =\{\varphi \in J^I \vert : \exists f \in B^A : (\varphi, f)\in \mathbf{bpol}_\bullet(\bfa, \bfb) \}$. We denote by $[\varphi]_b$ the morphism corresponding to $\varphi$
\pg (iii) The identity on  p-site $\bfa\equiv \ppar{I, A, (A_i)}$  in $\mathbf{bpol}_{\bullet/b}$ is the identity $id_I$.
\pg
It is worth noting that in this category an isomorphism $\bfa \to \bfb$ exists if and only if $\bfa$ and $\bfb$ are B-isomorphic (Definition \ref{BSisomorphism} ). 

\pg
The functor $\scrK: \mathbf{bpol}_\bullet \to \mathbf{P}^2_\bullet$  descends to  a functor denoted also $\scrK: \mathbf{bpol}_{\bullet/b} \to \mathbf{P}^2_\bullet$, by positing on object $\bfa$, $\scrK(\bfa) := (I, \kappa_\bfa)$, and on a morphism $[\varphi] \in \mathbf{bpol}_{\bullet/b} (\bfa, \bfb)$, $\scrK ([\varphi]) : = \varphi$. 

\pg The functor $\scrCa: \mathbf{P}^2_\bullet \to \mathbf{bpol}_\bullet$ descends (by taking quotient) to a functor  $\scrCa_b: \mathbf{P}^2_\bullet \to \mathbf{bpol}_{\bullet/b}$, by positing $\scrCa_b(I, \calX) = \bca (I, \calX)$   and on  a morphism $\varphi\in \mathbf{P}^2_\bullet ((I, \calX), (J, \calY))$, $\scrCa_b(\varphi) :=[\varphi]_b$ (in view of Proposition \ref{BSexistence2}, $\varphi \in \mathbf{BPOL}_\bullet(\bca(I, \calX), \bca(J, \calY)$). The verification of the  functor properties  of $\scrCa_b$ and $\scrK$ is easy routine.

 \pg The  following statement  is the main structure theorem concerning  the B-political theory:
 \begin{thm}\label{thm2B}
$\begin{tikzcd}
{\mathbf{bpol}_{\bullet/b}}\arrow[r, "{\scrK} ", yshift =0.7ex]&{\mathbf{P}^2_\bullet}\arrow[l, " \scrCa_b"  , yshift=-0.7ex]
\end{tikzcd}$
is an equivalence of categories
\end{thm}
\emph{Proof.}  Since  ${\scrK} \circ \scrCa_b =\mathrm{1}_{\mathrm{P}^2_\bullet} $,  it remains to prove that there is a natural isomorphism  between  the two functors $\scrCa_b \circ {\scrK}$   and $\mathrm{1}_{\mathrm{BPOL}}$, in symbols : $ \scrCa_b \circ {\scrK}\cong \mathrm{1}_{\mathbf{bpol}_{\bullet/b}}$. 
 The following diagram 
\begin{equation}
\begin{tikzcd}[column sep=large]
{\bfa}\arrow[d, "{\pi}_\bfa" ']\arrow[r, "{\varphi, f}" ] &{\bfb}\arrow[d, "{\pi}_\bfb "] \\
{\scrCa_b\circ \scrK(\bfa)}\arrow[r,  " {\scrCa_b \circ \scrK (\varphi, f)}" ']&{\scrCa_b\circ \scrK(\bfb)}
\end{tikzcd}
\end{equation}
is commutative since the lower morphism $\scrCa_b\circ \scrK (\varphi, f)$  is equal to $[\varphi]$.
Therefore the transformation $ \pi =( \pi_\bfa)$  from $\mathrm{1}_{\mathbf{bpol}_{\bullet/b} }$ to $ \scrCa_b\circ \scrK$ is a natural transformation.  In fact  $\pi =(\pi_\bfa)$  is  a natural isomorphism:
indeed $\bfa$ and $\scrCa_b \circ \scrK (\bfa) \equiv \bca \circ \kappa (\bfa)$ are BG-isomorphic (Definition \ref{Bisomorphism}), since they have the same base, and their Knit  are both equal to $\kappa_\bfa$, in particular  they are B-isomorphic (Definition \ref{BSisomorphism}), hence they are isomorphic in  $\mathbf{bpol}_{\bullet/b}$.   We conclude that  $\mathbf{bpol}_{\bullet/b} \simeq \mathbf{P}^2_\bullet$. \hfill $\Box$
\paragraph{$\mathbf{BPOL}_{\bullet}$ and isomorphism with $\mathbf{P}^2_\bullet$} --
We introduce on the p-sites a partial preorder relation $ \to_B$ as follows: 
 $\bfa \to_B\bfa'$ if and only if  $\mathrm{Base}\,(\bfa)= \mathrm{Base}\,(\bfa')$  and if there exists  a BG-map $ f : \mathrm{Ground}(\,\bfa)\to \mathrm{Ground}(\, \bfa') $, and the associated  equivalence relation $\sim_B$  defined by $\bfa \sim_B \bfa'$ if and only if $\bfa \to_B \bfa'$ and $\bfa'\to_B \bfa$.  Clearly $\bfa\sim_B\bfa'$ if and only if  $\bfa$ and $\bfa'$ are BG-isomorphic (Definition \ref{Bisomorphism}). 
 
\pg It is easily checked that for $\bfa \sim_B \bfa'$ and $\bfb\sim_B\bfb'$ one has $ \mathbf{bpol}_{\bullet/b}(\bfa, \bfb) =\mathbf{bpol}_{\bullet/b}(\bfa', \bfb')$

\pg We define a new category $\mathbf{BPOL}_\bullet$ the objects of which are equivalence classes of $\sim_B$ and where a morphism of codomain $\dot \bfa$ and codomain $\dot \bfb$ is any element of $\mathbf{bpol}_{\bullet /b}$, in symbols : $\mathbf{BPOL}_\bullet(\dot\bfa,\dot\bfb) = \mathbf{bpol}_{\bullet/b}(\bfa, \bfb)$. 
 
\pg The functor $\scrK$ descends to a functor (denoted with the same symbol)  $\scrK: \mathbf{BPOL}_\bullet\to  \mathbf{P}^2_\bullet$ and the functor $\scrCa_b$ descends (by taking quotient) to a functor $\scrCa_B :  \mathbf{P}^2_\bullet \to \mathbf{BPOL}_\bullet$.

\begin{thm}\label{thmBposet} The category $\mathbf{BPOL}_\bullet$ is  isomorphic to $\mathbf{P}^2_\bullet$
\end{thm}
\emph{Proof}.  $\scrK \circ \scrCa_B =\mathrm{1}_{\mathbf{P}^2_\bullet} $ and $\scrCa_B \circ \scrK =\mathrm{1}_{\mathbf{BPOL}_\bullet}$ \hfill $\Box$

\subsubsection{Category $\mathbf{spol}_\bullet$ }
The construction of this subsection is parallel to that of Subsection \ref{pol-B},  the B-concepts being replaced par the corresponding S-concepts. We start by describing some specific aspects of the S-case.

\paragraph{Preliminaries}
The \emph{effective ground} of a p-site $\bfa\equiv \ppar{I, A, (A_i)}$ is the set $A_\ast : =\cup_{i\in I} A_i = \{x\in A\vert\, \pi_\bfa(x)\neq\emptyset \}$. The \emph{effective p-site} associated with $\bfa$ is the p-site $\bfa_\ast:= \ppar{I, A_\ast, (A_i)}$. It is characterized by the parting function $\pi_{\bfa_\ast} = {\pi_\bfa}_{\vert {A_\ast}}$. Clearly $\kappa_{\bfa_\ast}=\kappa_\bfa\setminus\{\emptyset\}$ and $\sigma_{\bfa_{\ast}}= \sigma_\bfa$.
\pg
An S-map $(\varphi, f): \bfa\to \bfb$ induces a  by restriction on domain  an codomain and  function   $f_\ast : A_\ast \to B_\ast$ defined by $f_\ast (a)= f(a) (a\in A_\ast)$. $(\varphi, f_\ast)$ is clearly an S-map of domain $\bfa_\ast$ and codomain $\bfb_\ast$.  
Two SG-maps $(\varphi,f)$ and $(\varphi', f')$ with domain $\bfa$ and codomain $\bfb$ are  said to be {S-equivalent}  denoted  $(\varphi, f) \equiv_S (\varphi', f')$ if   $\varphi = \varphi'$ and $f_\ast= f'_\ast$. The equivalence class of $(\varphi,f)$ is denoted $[\varphi, f]_S$. 
\paragraph{Definition of the category $\mathbf{spol}_\bullet$} -- The objects of the category    $\mathbf{spol}_\bullet$ are  p-sites and the morphisms of domain $\bfa$ and codomain $\bfb$ in this category,  denoted $\mathbf{spol}_\bullet (\bfa, \bfb)$,  consist of the \emph{quotient  of the set of  all S-maps of domain $\bfa$ and codomain $\bfb$, with respect to the equivalence relation $\equiv_S$}.   

\pg The association $[\varphi, f]_S \mapsto (\varphi, f_\ast)$ is a bijection of
the set $\mathbf{spol}_\bullet(\bfa, \bfb)$ to the set of SG-maps of domain $\bfa_\ast$ and codomain $\bfb_\ast$. The two sets will be identified via this bijection.  Using this identification we have  $\mathbf{spol}_\bullet(\bfa, \bfb) = \mathbf{spol}_\bullet(\bfa_\ast, \bfb_\ast)$.
In particular the identity on $\bfa$  is the equivalence class  of  $id_A$.

\pg It easily seen that the category  $\mathbf{spol}_\bullet$ is well defined thanks to Proposition \ref{category1} and Fact \ref{category2}.

\paragraph{The Nerve functor $\scrS$}
The association $\scrS :\bfa\mapsto (I, \sigma_\bfa)$  for objects and $\scrS : (\varphi, f) \mapsto \varphi$ for arrows represented by the diagram : 

\begin{equation}
\begin{tikzcd}
 {\mathbf{SIM}_\bullet}\\
{\mathbf{spol}_\bullet }\arrow[u, "\scrS"]
\end{tikzcd}
\hskip12pt
\begin{tikzcd}[column sep=large]
{(I, \sigma_\bfa)}\arrow[r,  "{\varphi} " ]&{(J,\sigma_\bfb)}\\
{\bfa}\arrow[r, "{\varphi, f}" ]\arrow[u, maps to] &{\bfb}\arrow[u, maps to]
\end{tikzcd}
\end{equation}
is easily proven to be a functor of domain $\mathbf{spol}_\bullet$ and codomain $mathbf{SIM}_\bullet$  In the opposite direction, we shall denote  by the same symbol the restriction of the functor $\scrCa$ to $ {\mathbf{SIM}_\bullet}$, as shown in the following diagram:

\begin{equation}
\begin{tikzcd}
 {\mathbf{SIM}_\bullet}\arrow[d, "\scrCa"]\\
{\mathbf{spol}_\bullet }
\end{tikzcd}
\hskip12pt
\begin{tikzcd}[column sep=large]
{(I,\calX)} \arrow[d, maps to]\arrow[r,hook,   "{\varphi} " ]&{(J,\calY)}\arrow[d, maps to]\\
{\bca(I,\calX)}\arrow[r,  "{\varphi,{\hat\varphi}^\calX_\calY} "] &{ \bca(J,\calY)} 
\end{tikzcd}
\end{equation}

\paragraph{The morphism ${\overset{\ast}\pi_\bfa}$} -- The  parting function  of $\bfa$, $\pi_{\bfa}$, induces the function ${\overset{\ast}\pi_\bfa} : A_\ast \to \sigma_{\bfa}$ by positing  ${\overset{\ast}\pi_\bfa} (x) : = \pi_\bfa (x)$ $(x\in A_\ast)$. 

\begin{lem}\label{transformationS2}
 ${\overset{\ast}\pi_\bfa}$  is a BG - map of domain $\bfa_\ast$ and codomain $\bca(I, \sigma_\bfa)$.
\end{lem}
\emph{Proof.} Put $\bfb := \bca(I, \sigma_\bfa)$ and $f: {\overset{\ast}\pi_\bfa}$. Since $\pi_\bfb = i^{\sigma_\bfa}_{\calP(I)}$ the inclusion, one has for any $x\in A_\ast$,  $\pi_{\bfa_\ast} (x)= \pi_{\bfa}(x) =  i^{\sigma_\bfa}_{\calP(I)} (\pi_{\bfa}(x))= \pi_\bfb (f(x))$. We conclude that $\pi_{\bfa_\ast}= \pi_{\bfb} \circ f$, equivalently that $f$ is a BG-map. \hfill $\Box$

\begin{cor}  $(id_I, {\overset{\ast}\pi_\bfa}) \in \mathbf{spol}_\bullet\big(\bfa, \bca(I, \sigma_\bfa)\big)$
\end{cor} 
\pg We are going to prove that the family indexed by $\bfa$, $\overset{\ast}\pi \equiv \big((id_I,{\overset{\ast}\pi_\bfa})\big)$, is a natural transformation, precisely  a map  of domain the functor   $\mathrm{1}_{\mathbf{spol}_\bullet }$   and codomain the functor  $\scrCa\circ \scrS$.  

\begin{lem}\label{transformationS} $\overset{\ast}\pi$  is a natural transformation from $\mathrm{1}_{\mathbf{spol}_\bullet }$  to $ \scrCa\circ \scrS$ also written in symbols  $ \mathrm{1}_{\mathbf{spol}_\bullet}\overset{{\overset{\ast}{\pi}}} \Longrightarrow \scrCa\circ \scrS$ 
\end{lem}
\emph{Proof.} Since $\scrCa\circ \scrS (\varphi,f)= \bca\circ \sigma (\varphi,f) = (\varphi, \hat \varphi ^{\sigma_\bfa}_{\sigma_\bfb})$, the following diagram is commutative :
\begin{equation}
\begin{tikzcd}[column sep=large]
{\bfa}\arrow[d, "{\overset{\ast}\pi_\bfa}" ']\arrow[r, "{\varphi, f}" ] &{\bfb}\arrow[d, "{\overset{\ast} \pi_ \bfb} "] \\
{\scrCa\circ \scrS(\bfa)}\arrow[r,  "{ \scrCa\circ \scrS (\varphi, f)} " ']&{\scrCa\circ \scrS(\bfb)}
\end{tikzcd}
\end{equation}
\hfill $\Box$

\paragraph{$\mathbf{spol}_{\bullet/s}$ and equivalence with $\mathbf{SIM}_\bullet$}
- The category  $\mathbf{spol}_{\bullet/s}$ has p-sites as objects, and has  as morphisms of domain $\bfa$ and codomain $\bfb$ the set: 
\begin{equation}  
 =\{\varphi \in J^I \vert : \exists f \in B^A : (\varphi, f)\in \mathbf{spol}_\bullet(\bfa, \bfb) \}
\end{equation}
\pg  To any function $\varphi : I \to J$, one can associate $[\varphi]_s:=\{(\varphi, f)\vert\,  (\varphi, f) \in  \mathbf{spol}_\bullet(\bfa, \bfb)\}$.
 In view of \ref{BSexistence2}  $[\varphi]_s \neq \emptyset$ if and only if $ \hat \varphi (\sigma_\bfa) \subset \sigma_\bfb$.  It follows that 
$ \mathbf{spol}_{\bullet/s}(\bfa, \bfb) =\{\varphi \in J^I\vert \hat \varphi ( (\sigma_\bfa) \subset \sigma_\bfb\}$
 \pg
 We then define:
\pg - a functor  $\scrS : \mathbf{spol}_{\bullet/s} \to \mathbf{SIM}_\bullet$ by  positing $\scrS (\bfa) : = (I, \sigma_\bfa)$ on objects,  and 
$\scrS (\varphi) : =\varphi$ on morphisms
\pg -  a functor $\scrCa_s :  \mathbf{SIM}_\bullet\to \mathbf{spol}_{\bullet/s} $ by positing  $\scrCa_s (I, \calX) := \bca(I, \calX)$ on objects and $\scrCa_s(\varphi):= \varphi$ on morphisms

\pg  The  following statement  is the main structure theorem concerning  the S-political theory:

\begin{thm}
$\begin{tikzcd}
{\mathbf{spol}_{\bullet/s}}\arrow[r, "{\scrS} ", yshift =0.7ex]&{\mathbf{SIM}_\bullet}\arrow[l, "\scrCa_s"  , yshift=-0.7ex]
\end{tikzcd}$
is an equivalence of categories
\end{thm}
\emph{Proof.} The proof is similar to that of Theorem \ref{thm2B}, and where reference to Lemma \ref{transformationB2} is replaced by reference to Lemma \ref{transformationS2}. \hfill $\Box$
\paragraph{$\mathbf{SPOL}_\bullet$ and isomorphism with $\mathbf{SIM}_\bullet$} --
We  introduce  the equivalence relation  $\sim_S$ on the class of p-sites as follows: $\bfa \sim_S\bfa' $ if and only if $\bfa$ and $\bfa'$ are S-isomorphic (Definition \ref{BSisomorphism}). The definition of  $\mathbf{SPOL}_\bullet $  follows the same pattern as the that of  $\mathbf{BPOL}_\bullet $ 
by replacing $\sim_B$ by $\sim_S$. Finally we have:
\begin{thm}\label{thmSposet} The category $\mathbf{SPOL}_\bullet$ is  isomorphic to $\mathbf{SIM}_\bullet$
\end{thm}
\emph{Proof}. Same proof as Theorem \ref{Bposet1} where $\sim_B$ is replaced by $\sim_S$. \hfill $\Box$

\section{Scope of the Model and Political Applications}

\paragraph{The Categorical Framework}
-- The model proposed here establishes a categorical framework for the analysis of political intertwining and transformation. Working within this framework involves identifying two fundamental objects: (1) the site (situs, space, or place) in which political life unfolds, and (2) the general form of transition from one site to another. As is well known, any category can be regarded as a collection of morphisms governed by a principle of composition. In this perspective, the model does not directly represent empirical political dynamics; rather, it specifies the admissible transformations between political sites and the compositional rules ensuring their coherence.

\paragraph{Observables and Theory}

Underlying the two-level construction -- Configuration/Foundation -- is the idea that the configuration of the polity is, in part, observable, while the p-site at its foundation is a hypothetical construct posited to account for the observed structure. Since multiple p-sites may underlie a given p-formation, the political analyst must formulate a hypothesis concerning the ground and profile of the polity. This includes defining the ground and describing the aspiration sets of the parties. Some of this information may be publicly available through party declarations, while other aspects may be inferred from observation and empirical research.

\paragraph{ What the Model Is and Is Not}

The categorical construction provides the conceptual basis for a theory of political transformation but does not itself constitute a model of political action or behavior. This limitation reflects the deliberately abstract nature of its objects: a political party is defined solely as a set of desirable states (aspirations). It is not endowed with power, strategy, or agency, nor is it guided by preferences or rationality criteria, as in Game Theory, Social Choice Theory, or Decision Theory. This abstraction is intentional. By suspending empirical assumptions, it isolates the structural space within which behavioral models may subsequently be embedded. In this sense, the category serves as a \textit{conceptual arena} rather than a predictive mechanism. Example \ref{gallopolis2} illustrates how such an arena can nevertheless function as a paradigm for application.

Consider a case in which an event occurs within a political world represented by a p-site $\bfa$. The categorical framework requires only that, following this event, some morphism with domain $\bfa$ be selected. The framework itself does not prescribe how this selection occurs; that role belongs to a behavioral model, if one is specified. In the Gallopolis case, the implicit behavioral model stipulates that political bargaining results in the formation of a winning coalition  -- one capable of securing a parliamentary majority. Here, the bargaining process functions as the mechanism selecting a particular morphism of the category, namely, the one corresponding to a viable and victorious coalition.

It must be emphasized, however, that this implicit behavioral model depends on elements external to the category: (i) the attribution to each party of the weight corresponding to its parliamentary representation (its number of seats); (ii) the institutional rules governing parliamentary decision-making (as codified in the constitution); and (iii) the bargaining process itself, motivated -- at least in par -- by tactical advantage among negotiating actors. These elements belong not to the categorical structure but to the empirical domain of political practices and institutions. The category thus defines the \textit{conditions of possibility} for transformation, while the concrete realization of political behavior unfolds within the empirical dynamics of bargaining, strategy, and institutional constraint.

\section{Conclusion}

We  have  constructed a two-level  model of  political organisation  and transformation. At the Configuration level we defined the Category  $\mathbf{P}^2_\bullet$  and its subcategory   $\mathbf{SIM}_\bullet$. At the Foundation level we defined  two  categories that both have p-sites as objects,  namely $\mathbf{spol}_\bullet$ with   S-maps as morphisms, and $\mathbf{bpol}_\bullet$ with B-maps as morphisms, Two Functors  the Knit and the Nerve have been defined; they go from the P-foundation to the P-configuration.  The former determines the latter. The  Knit maps the p-sites to  p-formations, and  the morphisms of $\mathbf{bpol}_\bullet$ to morphisms of  $\mathbf{P}^2_\bullet$,  while the Nerve maps the p-sites  to p-structures and the morphisms  of $\mathbf{spol}_\bullet$ to morphisms of  $\mathbf{SIM}_\bullet$.  Our study establishes that the structure of $\mathbf{bpol}_\bullet$ is comparable to that of $\mathbf{P}^2_\bullet$, while the structure of $\mathbf{spol}_\bullet$ is comparable to that of $\mathbf{SIM}_\bullet$. 
\p
The morphisms  of  both  $\mathbf{spol}_\bullet$  and $\mathbf{bpol}_\bullet$ defined in this part of the research are built on functions between sets.  Such morphisms can represent only union  of parties in the Base  and merger of states on the Ground. It follows that, by the very nature of  our morphisms, the level of conflict in the codomain is lower than that of  the domain. This feature of our model in the present work  justifies the subtitle \emph{Hom\'onoia}. In  what follows, we present an empirical example that shows that, even in a highly regulated political entity, the transformation of the political structure may include scissions as well as unions. 
\p

 \paragraph{Concluding Example: The Macron Transition}

The election of Emmanuel Macron as President of the French Republic in May 2017 marks a key transition in French political dynamics. The event confirmed the persistence of the emergent force \textit{En Marche} while initiating a broader recomposition of the political landscape. Several centrist parties experienced internal divisions, and substantial fractions of their members joined the new formation.

In categorical terms, this episode corresponds to a transformation between two political sites, depicted on the left and right of Figure \ref{macron}. Both configurations, simplified to five members, are defined by their viability relations (the simplexes). The morphisms between them represent the admissible transformations associated with the transition.

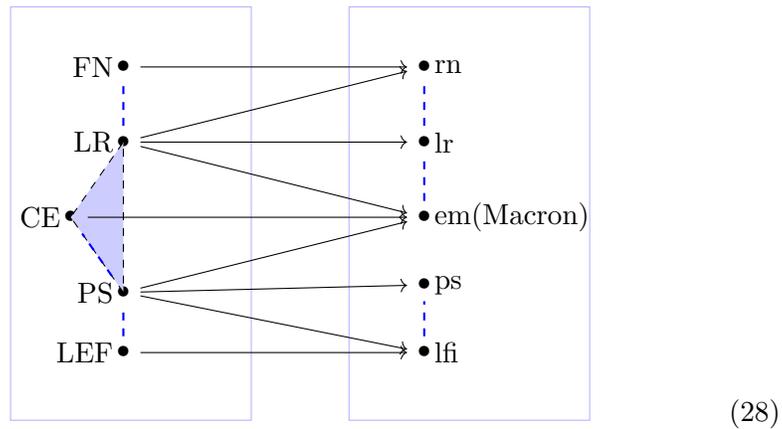
\begin{equation}\label{macron}
\begin{tikzpicture}[scale=1]
\draw[blue!30] (-0.5,-3) rectangle (2.7,2.5);
\draw (1,-2.1)node(1a){$\bullet$} node [left] {LEF};
\draw (1,-1.3)node(2a){$\bullet$} node [left] {PS};
\draw (0.3,-0.3)node(3a){$\bullet$} node [left] {CE};
\draw (1, +0.7)node(4a){$\bullet$} node [left] {LR};
\draw (1,+ 1.7)node(5a){$\bullet$} node [left] {FN};
\draw[thick,blue,dashed] (1a)--(2a);
\draw[thick,blue,dashed] (4a)--(5a);
\draw[fill=blue!20, dashed] (1,-1.3)--(0.3,-0.3)--(1, +0.7)--cycle;
\draw[thick,blue,dashed] (2a)--(3a);
\draw[blue!30] (4,-3) rectangle (7.2,2.5);
 \draw (5,-2.1)node(1b){$\bullet$} node [right] {lfi};
\draw (5,-1.2)node(2b){$\bullet$} node [right]{ps};
\draw (5,-0.3)node(3b){$\bullet$} node [right]{em \!\!(Macron)};
\draw (5,+0.7)node(4b){$\bullet$} node [right]{lr};
\draw (5,+ 1.7)node(5b){$\bullet$} node [right] {rn};
\draw [->](1a) edge (1b);
\draw [->](2a) edge (1b);
\draw [->](2a) edge (2b);
\draw [->](2a) edge (3b);
\draw [->](3a) edge (3b);
\draw [->](4a) edge (3b);
\draw [->](4a) edge (4b);
\draw [->](4a) edge (5b);
\draw [->](5a) edge (5b);
\draw[thick, blue, dashed] (1b)--(2b);
\draw [thick, blue, dashed] (3b)--(4b)--(5b);
\end{tikzpicture}
\end{equation}

Because a splitting party contributes to multiple successor entities, the mapping between the two configurations is inherently \emph{multivalued}. Capturing this process requires extending the categorical framework to include morphisms defined by multivalued functions on both the base and the ground. Such an extension generalizes the model to encompass transformations involving division, merger, or realignment -- core features of any political evolution.

\newpage
\appendix
\label{appendixa}
\begin{center}
{\huge \bf Appendix}
\end{center}
\section{ 
Delegation in the categorical framework}\label{delegation}

The notion of friendly delegation in the context of a political structure has been introduced  in \cite{AK2019}.  In that model, a political structure is represented by a simplicial complex, and  a  friendly delegation is a  move by some member that  consists in  delegating his role to another agent politically more central than himself. Such a move drives the political structure toward a less conflictual one.   In the following section we show how the notion of delegation can be viewed as a map  of the Configuration category.

\subsection{Delegation as a map of a P-configuration}\label{delegation1}
  
 \p
Let $I$ be a finite set, interpreted as a Base of some polity, and let $i_0, j_0 \in I$, $i_0\neq j_0$ be two distinct members of the base $I$. The \emph{delegation} from $i_0$ to $j_0$ is the function $\delta^{i_0}_{j_0} : I \to I $ defined by $\delta(i_0)=j_0$ and  $\delta(i)=i$ for all $i\neq i_0$. $i_0$ is the \emph{delegating} agent  $j_0$ is the \emph{delegate}. Let $(I,\calE)$ be a p-structure (Definition \ref{defpformation},  
the delegation $\delta^{i_0}_{j_0}$  in $(I,\calE)$ is said to be : 
\pg - \emph{simplicial}  in $(I, \calE)$ iff : $\!\!\forall \!s\in \!\calE$,  $i_0\in \!s \!\Rightarrow \{j_0\}\! \setminus\!\{i_0\} \!\cup s\in \calE$, equivalently : $\!\!\forall s\in \calE$, $\hat\delta(s)\in \calE$  
\pg-  \emph{friendly}  in $(I, \calE)$  iff : $\!\!\forall \!s\in\! \calE$,  $i_0\in s \Rightarrow \{j_0\} \cup s\in \calE$,  equivalently :
 $\forall s\in \calE$, $\hat\delta(s)\cup s\in \calE$  
  
\pg  In the terminology of discrete topology, a delegation $\delta$ is simplicial  in $(I,\calE)$ if and only if  $\delta : (I,\calE)\to (I,\calE)$ is a simplicial map, and it is 
 is friendly if and only if  $\delta: (I,\calE)\to (I,\calE)$ is contiguous to the identity  on $I$ (see \cite{Spanier}, or \cite{AK2019}). 
 In the terminology of the P-configuration category, a delegation $(\delta, \calE)$ is simplicial if and only if the function $\delta$ induces an SC-map from the p-structure $(I, \calE)$ into itself.

\begin{nota} $\calE^{-i_0} := \{s\in \calE\vert\,  i_0\notin s\}$ 
\end{nota}

\pg In \cite{AK2019}, the effectuation of a (friendly) delegation $\delta^{i_0}_{j_0}$ was understood as the removal of party $i_0$ and the passage from the original political structure  $(I,\calE)$ to  the new one $(I, \calE^{-i_0})$.  In the framework  of the P-configuration,  the effectuation of the delegation $\delta$ can be viewed as an action that operates on the p-formation $(I,\calE)$ and results in the p-formation $(I,\hat\delta\calE)$. 

\begin{fact} \label{fact1a}  For any delegation function  $\delta\equiv\delta^{i_0}_{j_0}:I\to I$, and any p-structure  $(I,\calE)$ one has 
$\hat\delta \calE^{-i_0}= \calE^{-i_0}$. 
\end{fact}

\pg What characterizes an simplicial delegation is asserted in the following:

\begin{pro} For any delegation function  $\delta\equiv\delta^{i_0}_{j_0}:I\to I$, and any p-structure  $(I,\calE)$, the following statements are equivalent :
\pg (i) $\delta$ is simplicial in $(I, \calE)$, 
\pg (ii)  $\delta : (I,\calE) \to (I, \calE)$ is SC-map,
\pg (iii) $\delta : (I,\calE) \to (I, \calE^{-i_0})$ is a BC-map,
\pg (iv) $\hat\delta(\calE)= \calE^{-i_0}$
\end{pro}


\begin{fact}\label{fact1b}
Any friendly delegation is simplicial.
\end{fact}

\pg Note that  the converse of Fact \ref{fact1b} is not true, as shown by the following:
\begin{exa}\label{notsofriendly}\rm  $I=\{1,2,3\}, \calE =\{12,13, 23, 1,2,3\}$ (boundary of a triangle), $\delta\equiv\delta^1_2 $. The delegation $\delta$ is not friendly since $13\in \calE$ and $123\notin \calE$.  However $\hat \delta \calE= \{23, 2,3\} =\calE^{-1}$, meaning that $\delta$ is simplicial. \hfill $\bigcirc$
\end{exa}

\pg In order to state what is specific about friendly delegations we first introduce the following:
\begin{nota}\rm   For any p-structure $(I,\calE)$, any  $i_0, j_0 \in I$, $i_0 \neq j_0$:   we put : $\calE^{i_0}_{j_0}:= \{s\in \calE\vert\, i_0\in s\Rightarrow j_0\in s\}$
\end{nota}

\begin{pro} \label{delegationchar} 
For any delegation function  $\delta\equiv\delta^{i_0}_{j_0}:I\to I$, and any p-structure  $(I,\calE)$, the following statements are equivalent :

\pg(i)The delegation  $\delta^{i_0}_{j_0}$ in $(I, \calE)$ is friendly, 
\pg(ii) $\calE=\calE^{i_0 \downarrow}_{j_0}$
\pg(iii) $\calE^{max}\subset \calE^{i_0}_{j_0}$
\end{pro}
\emph{Proof}. $(i)\Rightarrow(ii)$.  Let $s\in \calE$. If $i_0\notin s$ then $s\in \calE^{i_0}_{j_0}$. If $i_0\in s$, then  by assumption ${j_0}\cup s\in \calE$ and therefore  ${j_0}\cup s\in \calE^{i_0}_{j_0}$. It follows that $s\in \calE^{i_0 \downarrow}_{j_0}$. This proves  $\calE\subset \calE^{i_0 \downarrow}_{j_0}$. The opposite inclusion being always true, we proved  $\calE=\calE^{i_0 \downarrow}_{j_0}$.

\pg $(ii) \Leftrightarrow (iii)$ Let $s\in \calE^{max}$ then there exists  $t\in  \calE^{i_0}_{j_0}$ such that $s \subset t$. Since $s$ is maximal  for inclusion in $\calE$, $s=t$ and therefore $s\in \calE^{i_0}_{j_0}$. Conversely if  $\calE^{max}\subset \calE^{i_0}_{j_0}$ then $\calE^{i_0\downarrow}_{j_0} \supset \calE^{max\downarrow}=\calE$.
\pg $(iii)\Rightarrow(i)$. Conversely assume that  $\calE= {\calE}^{{i_0} \downarrow}_{j_0}$. Let $s \in \calE $, $i_0\in s$. If $j_0\in s$ then $\{j_0\}\cup s = s\in \calE$. If $j_0\notin s$ then there exists $t\in \calE^{i_0 \downarrow}_{j_0}$ such that $s\subset t$. Since $i_0\in s \subset t \in \calE^{i_0\downarrow}_{j_0}$, it follows that $j_0\in t$ and consequently ${\{ j_0 \} \cup s}\subset \in \calE$.  Again $\{j_0\} \cup  s \in \calE$. This proves  that $\delta_{j_0}^{i_0}$ is friendly in $(I, \calE)$. \hfill $\Box$

\subsection{Foundation  of a delegation on an invariant ground}\label{delegation2}
   The question arises   about the interpretation  of  the notion of delegation  (Subsection \ref{delegation1}) on the  ground i.e., on the level of the P-foundation. 

We say that  a  p-structure $(I,\calE)$ is \emph{$\sigma$-founded} on the  p-site $\bfa$ of ground $A$,  if $\bfa\in \mathrm{PSITE}(I, A)$  and $\sigma_\bfa=\calE$.
In view of Corollary \ref{canonical3}, any p-formation $(I,\calE)$ where $\calE$ is a $I$-simplicial complex, is $\sigma$-founded on a p-site $\bfa$ with  ground $\calA$ where  $\calE^{max}\subset \calA\subset \calE$, namely   the canonical p-site associated to $(I,\calA)$, $\bca(I, \calA)$,  Indeed one has $\kappa_\bfa=\calA$ and consequently $\sigma_\bfa=\kappa_\bfa^\downarrow=\calE$. 
 
The rest of this paragraph is devoted to understand what is specific about the $\sigma$-foundation of  friendly delegations.  
In \cite{AK2019}, the effectuation of a delegation $\equiv\delta^{i_0}_{j_0}$ is interpreted as the passage from $\calE$ to  $\hat\delta^{i_0}_{j_0}\calE \subset \calP(I\setminus\{i_0\})$;  in particular, the delegating agent is put aside from the resulting p-formation. If the delegation is simplicial and a fortiori friendly, the resulting p-structure is $(I,\calE^{-i_0})$. The latter expresses the simple withdrawal  of $i_0$  from the initial  p-structure without any other change.  At the level of the p-foundation the simple withdrawal of agent $i_0$ from the p-site is expressed by  a new profile that attributes   the empty set as aspirations to the withdraw party, while the aspirations of other parties  are kept unchanged.    This motivates the following:

 \begin{defi} \rm  For any p-site $\bfa =\ppar{I, A, (A_i)}$,  the  p-site that results from \emph{withdrawal} of $i_0$ is  $\bfa^{-i_0} := \ppar{ I,A,(B_i)}$ where $B_{i_0} :=\emptyset, B_i:= A_i \, (i\in I\setminus\{i_0\})$. \hfill$\triangledown$
\end{defi}

\g Note that in general, within  the p-site  $\delta^{i_0}_{j_0} \bfa= \ppar{I, A, (B_i)}$, image of $\bfa$ by $\delta^{i_0}_{j_0}$ (Definition and notation \ref{P-image}),   agent $i_0$ is silenced ($B_{i_0}= \emptyset$) but his/her aspirations are delegated   
to $j_0$ ($B_{j_0} =A_{i_0}\cup A_{j_0}$). The following  lemma details what  is implied by a delegation that results in the withdrawal of the delegating agent.

\begin{lem}\label {SCdelta} Let $ \delta^{i_0}_{j_0}: I\to I$ be a delegation function and let $\bfa\equiv \ppar{I,A,(A_i)}\in \mathrm{PSITE}(I,A)$.  The following are equivalent: 
 \pg (i) $\delta^{i_0}_{j_0}: \bfa \to \bfa^{-i_0}$ is a BP-map
  \pg (ii) $\delta^{i_0}_{j_0}\bfa=\bfa^{-i_0}$
 \pg (iii) $A_{i_0}\subset A_{j_0}$,
\pg (iv) $\delta^{i_0}_{j_0}: \bfa \to \bfa$ is an SP-map
\end{lem}

\begin{nota}\rm
We denote by $ \mathrm{PSITE}^{i_0}_{j_0}(I,A)$ the set of all p-sites on $(I, A)$ such that $A_{i_0}\subset A_{j_0}$.
\end{nota}

\g The p-sites of  $ \mathrm{PSITE}^{i_0}_{j_0}(I,A)$ are those where the delegation $\delta^{i_0}_{j_0}$ results in the withdrawal of $i_0$, the delegation does not modify the aspirations of the delegate.

\begin{pro} Let $\delta^{i_0}_{j_0}:I\to I$ be a delegation function. The following are equivalent:
\pg (i)The delegation $\delta^{i_0}_{j_0}$ on  $(I, \calE)$ is friendly, 
\pg(ii) $\bca(I, \calE^{max}) \in  \mathrm{PSITE}^{i_0}_{j_0}(I, \calE^{max})$,
\pg (iii)  There is  a set $A$,  a p-site $\bfa  \in \mathrm{PSITE}^{i_0}_{j_0}(I,A)$  such that $\sigma_\bfa =\calE$. 
\end{pro}
\emph{Proof}.  $(i) \Rightarrow (ii)$.  Put $\calA := \calE^{max}$ and $\bfa:= \bca(I, \calA)$. Clearly $\calA^{\downarrow}= \calE$, and in view of Proposition \ref{delegationchar} $\calA\subset \calE^{i_0}_{j_0}$.  Now, since $\bfa$ is canonical one has  $\pi_\bfa(s) = s  (s\in \calA) $.  It follows that for any $s\in \calA$,   if $i_0 \in \pi_\bfa(s)$, then $j_0 \in \pi_\bfa(s)$, equivalently $A_{i_0}\subset A_{j_0}$. We conclude that  $\bca(I, \calA)\in\mathrm{PSITE}^{i_0}_{j_0}(I,\calA)$.
\pg $(ii)\Rightarrow (iii)$. Trivial.
\pg $(iii)\Rightarrow (i)$. Assume that  $\sigma_\bfa =\calE$   and  $A_{i_0}\subset A_{j_0}$. If $i_0\in s \in \calE$ then $A_{s \cup \{j_0\}}=  A_s \neq\emptyset$ so that $s \cup \{j_0\}\in \calE$.  It follows that $\delta$ is friendly.
 \hfill$\Box$

\begin{example} (Illustration of a case of friendly delegation) In the upper part of the following figure, a delegation $\delta\equiv \delta^4_2$  is defined on the Base $I\equiv \{1,2,3,4\}$ and is shown to transport the p-structure $(I,\calE)$  onto its sub-structure  $\calF\equiv\hat\delta\calE$. This transformation $\sigma$-founded on a ground of 5 elements   by a p-site  $\bfa$ of profile $(A_i),  i=1,2,3,4$ that has Nerve $\calE$. The p-site $\bfa$  is BP-mapped  by $\delta^4_2$ to the p-site  $\bfb$ with Nerve $\calF=\hat\delta^4_2\calE$.

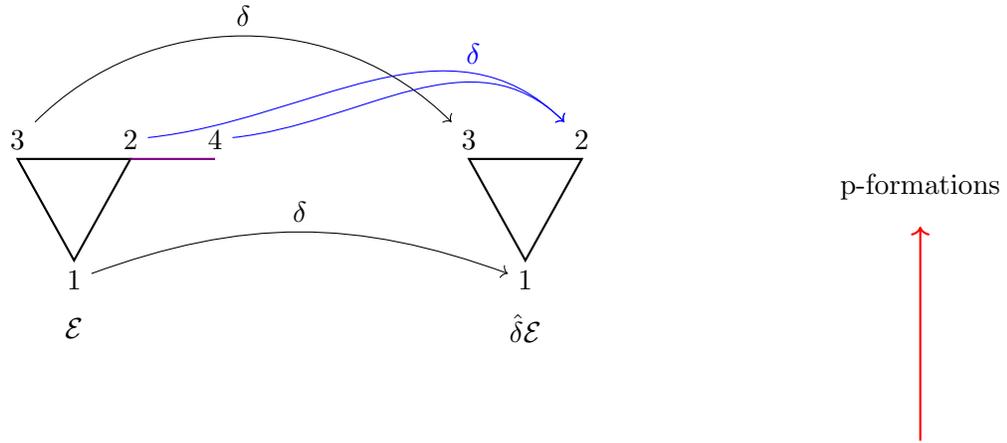
\begin{tikzpicture}[scale=0.75]
\draw[thick]  (-8,1) -- (-9,-0.8) -- (-10,1) -- cycle ;
\draw (-8,1) node (2a)[above] {$2$};
\draw[thick, violet] (-6.5,1) -- (-8,1);
\draw(-6.5,1) node (4a)[above]{$4$};
\draw (-9,-0.8) node(1a) [below] {$1$};
\draw (-10,1) node(3a) [above] {$3$};
\draw (-9,-2) node  {$\calE$};
\draw [thick] (0,1) -- (-1,-0.8) -- (-2,1) -- cycle ;
\draw (0,1) node(2b) [above] {$2$};
\draw (-1,-0.8) node(1b) [below] {$1$};
\draw (-2,1) node (3b) [above] {$3$};
\draw (-1,-2) node{$\hat\delta\calE$};
\draw [->](3a) edge [bend left=45] node[midway,above]{$\delta$}(3b);
\draw [->,blue](2a) to [out=6]node[near end,above]{$\delta$}(2b);
\draw [->,blue](4a) to [out=6](2b);
\draw [->] (1a) edge [bend left=20]node[midway,above]{$\delta$}(1b);

\node at (6,0.5){p-formations};
\draw[->, red, thick] (6,-4) -- (6,-0.2);
\end{tikzpicture}
\pp\pp

 
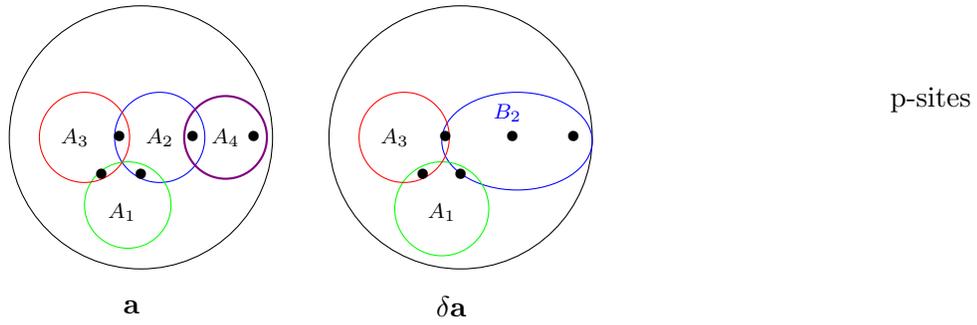
\begin{tikzpicture}[scale=0.25]
 \draw(2,0)circle[x radius= 7cm, y radius=7cm];
 \draw[blue] (3,0) circle (2.4);
 \draw[thick,violet](6.5,0)circle(2.2);
 \draw(6.5,0) node{$\scriptstyle  A_4$};
\draw[green] (1.3,-3.6) circle (2.3);
\draw[red] (-1,0) circle (2.4);
\draw(-1.5,0) node{$\scriptstyle  A_3$};
\draw(3,0) node{$\scriptstyle  A_2$};
\draw(1,-4) node{$ \scriptstyle A_1$};
\draw (-0.10,-2.0)node{$\bullet$}; 
\draw (0.85, 0)node{$\bullet$}; 
\draw (4.75, 0)node{$\bullet$}; 
\draw (2,-2)node{$\bullet$}; 
\draw (8,0)node{$\bullet$}; 
\draw (1.5,-9)node{$\bfa$};

\draw(19,0)circle  [x radius= 7cm, y radius=7cm];
 \draw[blue] (22,-0.2) circle [x radius=4.0cm, y radius=2.6cm];
\draw[green] (18.0,-3.8) circle (2.5);
\draw[red] (16,0) circle (2.4);
\draw(15.5,0) node{$\scriptstyle  A_3$};
\draw[blue](21.5,1.3) node{$\scriptstyle B_2$};
\draw(18,-4) node{$ \scriptstyle A_1$};
\node at (44,2){p-sites};

\draw (17,-2.0)node{$\bullet$}; 
\draw (18.2, 0)node{$\bullet$}; 
\draw (21.75, 0)node{$\bullet$}; 
\draw (19,-2)node{$\bullet$}; 
\draw (25,0)node{$\bullet$}; 
\draw (18.5,-9)node{$\delta\bfa$};

\end{tikzpicture}

\p\pg
Note that $\calE^{max}$ is composed of all the simplexes of dimension 1 of $\calE$, namely $\{12, 23,31,24\}$, and  it satisfies the appropriate condition $\calE^{max} \subset \calE^{4}_{2}$ that characterizes a friendly delegation. Therefore  the canonical p-site $\bca (I, \calE^{max})$ can be written  $\ppar {I, \calE^{max}, G_1, G_2,G_3, G_4}$ where $G_1=\{12,13\}$, $ G_2=\{23,12, 24\}, G_3=\{13,23\}$, $G_4=\{24\}$, so that in particular we have the inclusion $G_4\subset G_2$ as it should be.

\hfill $\bigcirc$

\end{example}

\end{document}